\documentclass[final]{siamart251216}
\usepackage{booktabs}
\usepackage{multirow}
\usepackage{amsmath}
\usepackage{amsfonts} 
\usepackage{braket,amssymb}
\usepackage{algorithmic}
\usepackage{array}
\usepackage{comment}
\usepackage{mathtools}
\usepackage{graphicx} %
\usepackage{subfig}
\usepackage{hyphenat}
\usepackage[sort]{cite}
\usepackage{anyfontsize}
\usepackage{hyperref}

\usepackage{tikz}
\usepackage{xcolor}
\usepackage{pgfplots}
\usepackage{pgfplotstable}
\pgfplotsset{compat=1.14}

\usepackage{pgfplots}
\theoremstyle{plain}

\pgfplotsset{compat=1.14}
\newsiamthm{assumption}{Assumption}

\title{A fast tangential rational Krylov subspace method with randomized sketching for large Lyapunov matrix equations.
\thanks{\textbf{Funding:} This work was supported in part by the National Science Foundation under grant DMS-2111496.}}
\author{
  Tianqi Zhang\thanks{School of Data Sciences, Zhejiang University of Finance and Economics, 
  18 Xueyuan Street, Hangzhou, Zhejiang 310018, China. Email: \texttt{tianqiz@zufe.edu.cn}} 
  \and
  Fei Xue\thanks{School of Mathematical and Statistical Sciences, Clemson University, 
  Clemson, SC 29634, USA. Email: \texttt{fxue@clemson.edu}}
}

\begin{document}

\maketitle

\begin{abstract}
We propose a tangential rational Krylov subspace method augmented with randomized sketching for solving continuous-time Lyapunov equations $AX + XA^T + BB^T = 0$, where $A$ is large and sparse and $B$ has a rank much higher than one. Numerical solutions require high-dimensional rational Krylov subspaces (e.g., of dimensions up to $\sim 1000$), leading to substantial computational cost and long runtime. Building on the adaptive strategies of [SIAM\,J.\,Matrix\,Anal.\,Appl., 35\,(2014), pp.\,476--498], which expand the subspace using the most relevant directions, we introduce several enhancements: selective orthogonalization of basis vectors and fast residual-norm estimation enabled by randomized sketching, and the option to use either adaptive distinct poles or a small set of repeated poles. We derive a new computable error bound for the approximate solutions and show that sketching does not degrade the convergence rate relative to the standard (non-sketched) algorithm. Numerical experiments demonstrate improved robustness and efficiency, with reduced runtimes and often higher accuracy than the standard tangential method.
\end{abstract}



\section{Introduction} This study concerns the numerical solution of the continuous-time Lyapunov equation of the form 
\begin{eqnarray}\label{lyapeqn}
    AX+XA^T+BB^T=0,
\end{eqnarray}
where $A\in \mathbb{R}^{n\times n}$ is a large sparse matrix, and $B \in \mathbb{R}^{n \times p}$ {has a rank $p \gg 1$ but $p \ll n$}. Lyapunov equations are critical in model reduction, stability and control theory. We assume, as is typical, that $A$ is stable, with all eigenvalues located strictly in the left half of the complex plane. This equation emerges from the linear system of first-order differential equations  $x'(t)=Ax(t)+B_1u(t)$ and $y(t)=B_2^Tx(t)$, where $B_1 \in \mathbb{R}^{n \times p_1}$, $B_2 \in \mathbb{R}^{n \times p_2}$, $x(t) \in \mathbb{R}^n$ represents the system state, $u(t) \in \mathbb{R}^{p_1}$ is the input, and $y(t) \in \mathbb{R}^{p_2}$ is the output. Solutions to the Lyapunov equations $AP+PA^T+B_1B_1^T=0$ and $A^TQ+QA+B_2B_2^T=0$ yield the controllability and observability Gramians, respectively \cite{gajic1995lyapunov}. The Gramians are critical for quantifying energy transfers within the system and play a pivotal role in model order reduction (MOR). By capturing the essential dynamical properties of the original large-scale system, MOR enables the construction of a significantly smaller, yet representative, model that is easier to explore numerically; see, e.g.,  \cite{Benner.Damm.2011,Baur.etal.2014,Benner.etal.2015} and references therein.
 
The numerical solution of large-scale Lyapunov equations has been a subject of research for over two decades, with a comprehensive overview in \cite{Simoncini2016}. Rational projection algorithms, such as the early extended Krylov subspace method (EKSM) \cite{Simoncini.2007} and the more advanced adaptive rational Krylov subspace method (RKSM) \cite{Druskin.Simoncini.2011}, are known to achieve significantly faster convergence compared to polynomial Krylov subspace methods \cite{Saad.1990,Jaimoukha.Kasenally.1994}. Among these, RKSM surpasses EKSM and the Alternating-Direction-Implicit (ADI) method \cite{Ellner.Wachspress.1986,Li.White.2002} in terms of convergence speed and robustness \cite{Druskin.etal.2011,Simoncini2016}. When the numbers of inputs and outputs, i.e., $p_1$ and $p_2$, are small (e.g., less than 5), block variants of EKSM or RKSM remain competitive. However, when $\max\{p_1,p_2\} \gg 1$, these block variants become computationally expensive due to the increased dimensionality of the required subspaces and the substantial cost needed to maintain an orthonormal basis of the subspaces.

To prevent the rapid expansion of the rational Krylov subspace dimension, a tangential variant of the adaptive rational Krylov subspace method (RKSM) has been introduced \cite{Druskin.etal.2014}. In each iteration, this algorithm identifies a unique new pole and the corresponding most significant vectors to augment the subspace efficiently. Experimental results show that, for specific problems, this approach yields an approximate solution with accuracy comparable to that of the standard block variant, while requiring a subspace of significantly smaller dimension. Therefore, it has been a powerful tool for solving Lyapunov equations where $\mathrm{rank}(B)\gg1$. Nevertheless, the computational cost becomes rather high when the subspace dimension reaches a few hundred, driven by the expense of orthogonalizing the basis vectors $V_m$ and the block $[AB,V_m]$ for pole selection and residual norm estimation. Also, this method relies on the basis vectors of the form $\{B,(A-s_1I)^{-1}Bd_1,\ldots,(A-s_kI)^{-1}Bd_k\}$, with \emph{distinct} poles $s_1,\ldots,s_k$. Therefore, if any two poles are nearly identical, the basis vectors become ill-conditioned and may undermine the accuracy of the computed solutions.

In this study, we introduce a new variant of the tangential algorithm that incorporates randomized sketching to enhance the efficiency and robustness of the existing method. Randomized sketching \cite{Woodruff.2014,Tropp.etal.2017,Martinsson.Tropp.2020} has gained traction in numerical linear algebra in recent years, valued for its ability to approximate solutions to linear least squares problems at reduced computational cost \cite{Sarlos.2006} and to support the use of non-orthogonal basis vectors in Krylov subspaces for iterative solutions of linear systems \cite{Nakatsukasa.Tropp.2024,BGS23}, eigenvalue problems \cite{Nakatsukasa.Tropp.2024}, {matrix Sylvester equations \cite{PalittaEtal25}}, and matrix function times a vector \cite{GS23,BG24,Cortinovis.etal.2024}. For rational Krylov subspace methods, randomized sketching has received limited attention because convergence is typically achieved in a moderate number of steps, with computational cost dominated by linear solves. However, randomized sketching proves particularly advantageous in our setting, where a high-dimensional subspace is required, to effectively lower the arithmetic cost. 

Our use of randomized sketching facilitates selective orthogonalization, fast computation of new poles, and efficient residual norm estimation. {Unlike the standard formulation used in existing TRKSM implementations, our method does not rely directly on basis vectors of the form \((A-s_iI)^{-1}Bd_i\) with distinct poles.} This flexibility allows the use of repeated poles, achieved by sparse LU factorizations, further boosting efficiency and robustness. In addition, we demonstrate that a key property of exactness for standard RKSM, associated with the compressed matrix $H_m = V_m^TAV_m$, extends to our new formulation, where the compressed matrix is $T_m = (SV_m)^\dagger SAV_m$ {($^\dagger$\,means the Moore-Penrose pseudoinverse, or left inverse here)} with the sketching matrix $S \in \mathbb{R}^{s \times n}$. The convergence analysis of RKSM for Lyapunov equations \cite{Druskin.etal.2011} rooted in the Faber-Dzhrbashyan (FD) rational functions \cite{Dzhrbashyan.1957},\cite[Ch.\,XIII, sect.\,3]{Suetin.1998} remains valid for the sketched variant. We provide deeper insights into the influence of poles on the errors of the approximate solutions.

We structure the remainder of this paper as follows. Section~2 provides a review of the standard tangential rational Krylov subspace method, randomized sketching, and selective orthogonalization techniques. In section~3, we introduce our novel sketched tangential RKSM, emphasizing its advancements in pole selection flexibility and efficient residual norm estimation. Section~4 establishes an exactness property for RKSM using sketched compressed matrices, demonstrating that the new variant exhibits convergence behavior consistent with the standard variant. Numerical experiments are given in section~5, and section~6 presents concluding remarks.

\section{Background}\label{sect:review} This section provides a quick review of the existing tangential RKSM, randomized sketching and selective orthogonalization. 

\subsection{Tangential RKSM} The tangential RKSM \cite{Druskin.etal.2014} for solving Lyapunov equations \eqref{lyapeqn} is as follows. Initially, it obtains the orthonormal basis vectors $V_p \in \mathbb{R}^{n \times p}$ from the reduced QR factorization of $B$, assuming that $\mathrm{rank}(B)=p$. Suppose that an orthonormal basis $V_{m_k}$ of dimension $m_k$ is available before the $k$-th step, generated from distinct poles $s_1,\ldots,s_{k-1} \in W \subset \mathbb{C}$. To find the new pole and the new vectors to expand the subspace, it defines the RKSM residual and the new pole as
\begin{eqnarray}\nonumber %
\qquad\,\,     {\mathcal R_k}(s) &=& (A-sI)V_{m_k}(H_{m_k}-sI)^{-1}V_{m_k}^TB-B \in \mathbb{R}^{n \times p}, \mbox{ with }H_{m_k}=V_{m_k}^TAV_{m_k}, \\ \label{maximizerrks}
    && \qquad\qquad\qquad \mbox{ and }s_k = \mathrm{argmax}_{s\in W}\|{\mathcal R_k}(s)\|,
\end{eqnarray}
where $W \subset \mathbb{C}$ is the domain for the poles of RKSM. The \emph{tangential vectors} $d_k \in \mathbb{R}^{p \times \ell_k}$ are defined as the leading $\ell_k$ ($1 \leq \ell_k \leq p$) right singular vectors of the $n\times p$ matrix ${\mathcal R_k}(s_k)$ associated with the largest singular values. The current subspace is expanded by orthogonalizing the $\ell_k$ vectors $(A-s_kI)^{-1}Bd_k$ against $V_{m_k}$ and orthonormalizing the updated vectors to obtain the expanded orthonormal basis $V_{m_{k+1}}=V_{m_k+\ell_k}$. %

The idea is to obtain the maximizer of $\|{\mathcal R_k}(s)\|$ as the new pole to expand the subspace so that $\max_{s\in W}\|{\mathcal R_{k+1}}(s)\| < \max_{s\in W}\|{\mathcal R_k}(s)\|$ holds in general. If $\|{\mathcal R_{k}}(s)\|$ is sufficiently small uniformly for all $s \in W$, then $(A-sI)V_{m_k}(H_{m_k}-sI)^{-1}V_{m_k}^TB \approx B$, i.e., $(A-sI)^{-1}B \approx V_{m_k}(H_{m_k}-sI)^{-1}V_{m_k}^TB$. Then any rational function $q_{j}(z) = c_0 + \sum_{i=1}^{j} c_i (z-\nu_i )^{-1}$ of $(j,j)$-type with distinct poles $\nu_i \in W$ satisfies the relation $$q_{j}(A)B  \approx \textstyle V_{m_k}\big(c_0 I + \sum_{i=1}^{j} c_i (H_{m_k}-\nu_i I)^{-1}\big)V_{m_k}^TB = V_{m_k} q_{j}(H_{m_k})V_{m_k}^TB.$$In particular, if $p = \ell_i \equiv 1$ for all $1 \leq i \leq k-1$ so that $m_k = k$, and if $j = k-1$ and $\nu_i = s_i$ for $1\leq i \leq k-1$, then $q_{k-1}(A)B = V_kq_{k-1}(H_{k})\|B\|e_1$. This is known as the exactness property and can be used to analyze the convergence of RKSM.

The computation of the maximizer in \eqref{maximizerrks} is expensive in its original definition due to the size of ${\mathcal R_k}(s)$ ($n$-by-$p$) for a large matrix $A \in \mathbb{R}^{n \times n}$. To speed up the evaluation, \cite{Druskin.etal.2014} first shows that $\|{\mathcal R_k}(s)\| = \big\|(AV_{m_k}-V_{m_k}H_{m_k})(H_{m_k}-sI)^{-1}V_{m_k}^TB \big\|$, and it gives the following identity (see Propositions 4.1 \& 6.1 therein)
\begin{eqnarray}\label{avminusvh}
  && \qquad AV_{m_k}-V_{m_k}H_{m_k} = [AV_{m_k},V_{m_k}]\left[\begin{array}{c}I_{m_k} \\ -H_{m_k}\end{array}\right] =  [AB,V_{m_k}] \times \\ \nonumber
    && \underbrace{\Bigg[\left[\begin{array}{cc}I_p & 0 \\ 0 & V_{m_k}^TB [d_1,\ldots,d_{k-1}]\end{array} \right]R_{m_k}^{-1}+\left[\begin{array}{c}0 \\ R_{m_k} S_{m_k} R_{m_k}^{-1}\end{array}\right],\,\left[\begin{array}{c}0 \\ I_{m_k}\end{array}\right]\Bigg]}_\text{$G$} \left[\begin{array}{c}I_{m_k} \\ -H_{m_k}\end{array}\right],
\end{eqnarray}
where $[d_1,\ldots,d_{k-1}] \in \mathbb{R}^{p \times (m_k-p)}$, $R_{m_k} \in \mathbb{R}^{m_k \times m_k}$ is the upper triangular factor of the reduced QR factorization of $\big[B,(A-s_1I)^{-1}B,\ldots,(A-s_{k-1}I)^{-1}B\big] \in \mathbb{R}^{n \times m_k}$, and $S_{m_k} = \mathrm{diag}(\underbrace{0,\ldots,0}_\text{$p$ times},\underbrace{s_1,\ldots,s_1}_\text{$\ell_1$ times},\ldots,\underbrace{s_{k-1},\ldots,s_{k-1}}_\text{$\ell_{k-1}$ times})$, so that $G \in \mathbb{R}^{(p+m_k)\times 2m_k}$. Let $\widehat{Q}_{p+m_k}\widehat{R}_{p+m_k} = [AB, V_{m_k}]$ and $Q_BR_B = B$ be reduced QR factorizations. Then 
\begin{eqnarray}\label{newrksnorm}
    \|{\mathcal R_k}(s)\| &=& \left\|[AB,V_{m_k}]G \left[\begin{array}{c}I_{m_k} \\ -H_{m_k}\end{array}\right] (H_{m_k}-sI)^{-1}V_{m_k}^TB \right\| \\ \nonumber
    &=& \Bigg\|\underbrace{\widehat{R}_{p+m_k}G \left[\begin{array}{c}I_{m_k} \\ -H_{m_k}\end{array}\right]}_\text{$\widehat{G} \in \mathbb{R}^{(p+m_k)\times m_k}$}(H_{m_k}-sI)^{-1}\underbrace{\left[\begin{array}{c}R_B \\ 0\end{array}\right]}_\text{$m_k \times p$}\Bigg\|.
\end{eqnarray}

In addition, the above computation and notation can be used to obtain the residual norm to implement the stopping criterion. Specifically, let $H_{m_k}Y_{m_k}+Y_{m_k}H_{m_k}^T+(V_{m_k}^TB)(V_{m_k}^TB)^T=0$ be the projected Lyapunov equation, and the residual norm of the original equation \eqref{lyapeqn} can be evaluated as
\begin{eqnarray}\label{resnorm}
    &&\left\|AV_{m_k}Y_{m_k}V_{m_k}^T + V_{m_k}Y_{m_k}V_{m_k}^TA^T+BB^T\right\| \\ \nonumber
    &=& \left\|\big(\widehat{R}_{p+m_k}G\big) \left[\begin{array}{cc} 0 & Y_{m_k} \\ Y_{m_k} & \big(V_{m_k}^TB\big)\big(V_{m_k}^TB\big)^T\end{array}\right]\big(\widehat{R}_{p+m_k}G\big)^T\right\|.
\end{eqnarray}

The bottleneck of the tangential RKSM is the storage and orthogonalization cost. {The derivation and the implementation of new pole selection and tangential vector computation \eqref{maximizerrks} and \eqref{newrksnorm}, and the equation residual norm evaluation \eqref{resnorm}, rely on the assumption that \emph{distinct} poles are adopted at each iteration step; see \cite[Propositions 4.2, 6.1, 6.2]{Druskin.etal.2014}.} The reduced QR factorizations of $\big[B,(A-s_1I)^{-1}Bd_1,\, \ldots,\, (A-s_{k-1}I)^{-1}Bd_{k-1}\big] \in \mathbb{R}^{n \times m_k}$ and $[AB,V_m] \in \mathbb{R}^{n \times (p+m_k)}$ are needed, typically performed by the modified Gram-Schmidt process with reorthogonalization. The two factorizations also need the storage of at least $p+2m_k$ vectors of $n$ entries before step~$k$. Therefore, this method becomes quite expensive when the subspace dimension $m_k$ increases to a few hundred. The matrix operations given in \eqref{avminusvh}, \eqref{newrksnorm}, and \eqref{resnorm} also require a non-trivial effort to implement. Also, these relations are derived on the basis of distinct poles, {which may lead to longer runtime than methods that reuse a small set of repeated poles through sparse LU factorizations.} In case two distinct poles are very close, the original basis vectors $\{B,(A-s_1I)^{-1}Bd_1,\ldots,(A-s_{k-1}I)^{-1}Bd_{k-1}\}$ are ill-conditioned, and the accuracy of approximate solutions could be negatively affected. As we shall see, these issues can be effectively mitigated by randomized sketching. 

\subsection{Randomized sketching} Randomized sketching {maps} high-dimensional data into a lower-dimensional space, maintaining an approximate isometry. Suppose that $V_m \in \mathbb{R}^{n \times m}$, where $\mathrm{rank}(V_m) = m \ll n$. {Given \(\epsilon\in(0,1)\), a fixed sketching matrix \(S\in\mathbb C^{s\times n}\) with \(m<s\ll n\) may be chosen, for example, as a dense random embedding with independent Gaussian or sign entries, a sparse dimension reduction map, or a structured embedding such as a subsampled randomized Fourier/trigonometric transform. With suitable \(s\), such matrices satisfy, with high probability,
\begin{eqnarray}\label{sketchassumption}
    (1-\epsilon)\|v\|_2^2 \leq \|Sv\|_2^2 \leq (1+\epsilon)\|v\|_2^2
\end{eqnarray}
holds uniformly for any $v \in \mathrm{span}(V_m)$.} The idea of such a mapping originates from the Johnson-Lindenstrauss lemma \cite{Johnson.Lindenstrauss.1984,Matousek.2008,Larsen.etal.2017}. The mapping $S$ is called an $\epsilon$-subspace embedding for $\mathrm{span}(V_m)$. In practice, we assume that \eqref{sketchassumption} is guaranteed. Examples of sketching include the subsampled random Fourier transform (SRFT) and the sparse dimension reduction map \cite{Martinsson.Tropp.2020}. The inequality \eqref{sketchassumption} can be extended to matrices. 

\begin{proposition}
    Given $F \in \mathbb{R}^{n \times m}$, let $S \in \mathbb{C}^{s \times n}$ be a corresponding randomized sketching matrix that satisfies \eqref{sketchassumption} for any vector $v \in \mathrm{span}(F)$. Then 
\begin{eqnarray}\label{F2normsineq}
    && (1-\epsilon)\|F\|_F^2 \leq \|SF\|_F^2 \leq (1+\epsilon)\|F\|_F^2, \\ \nonumber
    && (1-\epsilon)\|F\|_2^2 \leq \|SF\|_2^2 \leq (1+\epsilon)\|F\|_2^2
\end{eqnarray}
\end{proposition}
\begin{proof}
    Let $F_i$ denote the $i$-th column of $F$. From \eqref{sketchassumption}, we have
\begin{eqnarray}
    \textstyle (1-\epsilon)\|F\|_F^2 &=& \textstyle{(1-\epsilon) \sum_{i=1}^m \|F_i\|_2^2} \leq {\sum_{i=1}^m \|SF_i\|_2^2} = \|SF\|_F^2 \\ \nonumber
    & \leq & \textstyle {(1+\epsilon)\sum_{i=1}^m \|F_i\|_2^2}=(1+\epsilon)\|F\|_F^2.
\end{eqnarray}

Let $u_1$ and $\widehat{u}_1$ be the leading right singular vectors of $F$ and $SF$, respectively, so that $\|u_1\|_2 = \|\widehat{u}_1\|_2 = 1$, $\|Fu_1\|_2 = \|F\|_2$ and $\|SF\widehat{u}_1\| = \|SF\|_2$. Therefore, we have 
\begin{eqnarray}
    (1-\epsilon)\|F\|^2_2 &=& (1-\epsilon)\|Fu_1\|^2_2 \leq \|SFu_1\|^2_2 \leq \|SF\|_2^2 \\ \nonumber
     &=& \|SF\widehat{u}_1\|_2^2 \leq (1+\epsilon)\|F\widehat{u}_1\|_2^2 \leq (1+\epsilon)\|F\|_2^2. 
\end{eqnarray}
Both inequalities in \eqref{F2normsineq} are hence established.  
\end{proof}

Randomized sketching has gained increasing interest in matrix computations in recent years; see, e.g.,  \cite{Woodruff.2014,Martinsson.Tropp.2020}. With randomized sketching, several traditional methods that need to construct orthonormal basis vectors can proceed with a non-orthonormal yet not highly ill-conditioned basis, without significantly compromising the convergence rate or the accuracy of the approximate solutions \cite{Sarlos.2006,Nakatsukasa.Tropp.2024,BGS23,GS23,BG24}. %

To develop a non-orthonormal basis with randomized sketching at a reduced cost of orthogonalization, a particularly effective approach called selective orthogonalization has been proposed \cite{GS24}. Suppose that we want to orthogonalize vector $w \in \mathbb{R}^n$ against a set of linearly independent vectors $V_m\in \mathbb{R}^{n \times m}$ ($m < n$). The operation is achieved by $\widehat{w} = w - V_m V_m^\dagger w = (I-V_mV_m^\dagger)w$, where $(I-V_mV_m^\dagger)$ is an orthogonal projector onto $\mathrm{span}(V_m)^\perp$. With randomized sketching, we can similarly perform a full approximate orthogonalization by $\widehat{w} = w- V_m (SV_m)^\dagger (Sw)$ so that $S\widehat{w} \perp \mathrm{span}(SV_m)$. {In fact, as the sketching matrix induces a semidefinite inner product $\langle u, v \rangle_S = \langle Su, Sv \rangle$, this simply orthogonalizes $w$ against $V_m$ with respect to this inner product.} To reduce the cost, \emph{selective orthogonalization} first evaluates the vector $g = (SV_m)^\dagger (Sw) \in \mathbb{R}^m$, then it keeps the largest (in absolute value) $\ell$ elements of $g$ and zeros out the remaining $m-\ell$ elements to get the updated vector $\widehat{g}$ and obtain $\widehat{w} = w - V_m \widehat{g}$. This vector update involves only $\ell+1$ vectors instead of $m+1$ vectors. Since $SV_m \in \mathbb{R}^{s \times m}$, the evaluation of $(SV_m)^\dagger (Sw)$ is considerably less expensive than that of $V_m^\dagger w$ without sketching. With a proper choice of the \emph{depth of selective orthogonalization} $\ell$, this approach builds an expanding subspace with basis vectors whose condition number increases moderately with the subspace dimension growth. 

\section{A new tangential RKSM with randomized sketching}\label{sect:tgrksm} In this section, we propose a new randomly sketched tangential RKSM for solving large Lyapunov equations \eqref{lyapeqn} with $\mathrm{rank}(B) \gg 1$. As we shall see, the orthogonalization cost of the standard tangential RKSM can be reduced considerably by the sketched variant. 

\subsection{New pole selection and subspace expansion} As shown in section~\ref{sect:review}, tangential RKSM needs to select a new pole at each step by finding the maximizer of $\|{\mathcal R_k}(s)\|$ (see\,\eqref{maximizerrks}), achieved from a reduced QR factorization of $[AB,V_{m_k}] \in \mathbb{R}^{n \times (p+m_k)}$ at a rapidly increasing cost with the algorithm's progress. To streamline the maximization, suppose that we have a non-orthonormal basis $V_{m_k}$, and we can extend the definition of the RKSM residual accordingly as follows. 

\begin{proposition}\label{prop:res1}
    Given $V_{m_k} \in \mathbb{R}^{n \times m_k}$ of rank $m_k < n$, let $T_{m_k} = V_{m_k}^\dagger A V_{m_k}$ and ${\mathcal R_k}(s) = (A-sI)V_{m_k}(T_{m_k}-sI)^{-1}V_{m_k}^\dagger B-B$. If $\mathrm{span}(B) \subset \mathrm{span}(V_{m_k})$, then $${\mathcal R_k}(s) = (AV_{m_k}-V_{m_k} T_{m_k})(T_{m_k}-sI)^{-1}V_{m_k}^\dagger B \perp \mathrm{span}(V_{m_k})$$
\end{proposition}
\begin{proof}
    The proof is similar to that of \cite[Proposition 4.1]{Druskin.etal.2014}. We have
\begin{eqnarray}\nonumber
 {\mathcal R_k}(s) &=& AV_{m_k}(T_{m_k}-sI)^{-1}V_{m_k}^\dagger B- sV_{m_k}(T_{m_k}-sI)^{-1}V_{m_k}^\dagger B - B \\ \nonumber
 &=& AV_{m_k}(T_{m_k}-sI)^{-1}V_{m_k}^\dagger B + V_{m_k}(T_{m_k}-sI)(T_{m_k}-sI)^{-1}V_{m_k}^\dagger B \\ \nonumber
 && \quad - V_{m_k}T_{m_k}(T_{m_k}-sI)^{-1}V_{m_k}^\dagger B - B
 \end{eqnarray}
 To continue, we have
 \begin{eqnarray}
 {\mathcal R_k}(s) &=& \left(AV_{m_k}-V_{m_k}T_{m_k}\right)(T_{m_k}-sI)^{-1}V_{m_k}^\dagger B - \left(I-V_{m_k}V_{m_k}^\dagger\right)B \\ \nonumber
 &=& \left(AV_{m_k}-V_{m_k}T_{m_k}\right)(T_{m_k}-sI)^{-1}V_{m_k}^\dagger B.
\end{eqnarray}
The orthogonality follows from $V_{m_k}^\dagger (AV_{m_k}-V_{m_k}T_{m_k}) = 0$ as $T_{m_k} = V_{m_k}^\dagger AV_{m_k}$. 
\end{proof}

The RKSM residual ${\mathcal R_k}(s)$ defined above does not have an advantage in computational cost, because a reduced QR factorization of $V_{m_k}$ is still needed to realize the action of $V_{m_k}^\dagger$. To reduce the cost with randomized sketching, we can define the \emph{sketched RKSM residual} $S{\mathcal R_k}(s)$ as follows. 

\begin{proposition}\label{prop:res2}
    Given $V_{m_k} \in \mathbb{R}^{n \times m_k}$ of rank $m_k < n$, let $S \in \mathbb{C}^{s \times n}$ be a sketching matrix and assume that $SV_{m_k}$ is of full rank. Define {$T_{m_k} = (SV_{m_k})^\dagger SA V_{m_k}$} and $S{\mathcal R_k}(s) = S(A-sI)V_{m_k}(T_{m_k}-sI)^{-1}(SV_{m_k})^\dagger SB-SB$. If $\mathrm{span}(B) \subset \mathrm{span}(V_{m_k})$ $$S{\mathcal R_k}(s) = S(AV_{m_k}-V_{m_k} T_{m_k})(T_{m_k}-sI)^{-1}(SV_{m_k})^\dagger SB \perp \mathrm{span} (SV_{m_k})$$
\end{proposition}

The proof is omitted here since it is almost identical to that of Proposition \ref{prop:res1}. 

The evaluation of $S{\mathcal R_k}(s)$ is considerably cheaper than that of ${\mathcal R_k}(s)$ defined in Proposition \ref{prop:res1}. The matrix $T_{m_k} = (SV_{m_k})^\dagger (SAV_{m_k})$ is evaluated by a reduced QR factorization of $SV_{m_k}$, which is accumulated through the $k-1$ steps already taken. The matrix $SAV_{m_k}-SV_{m_k}T_{m_k} \in \mathbb{R}^{s \times m_k}$ can be evaluated efficiently by keeping $SAV_{m_k}, SV_{m_k} \in \mathbb{R}^{s \times m_k}$ during the $k-1$ steps taken, without storing $AV_{m_k}$ or any $n \times m_k$ matrices. The evaluation of $(SV_{m_k})^\dagger SB$ is also easy. Since we usually have $B = V_p B_p$ for some $B_p \in \mathbb{R}^{p \times p}$, $(SV_{m_k})^\dagger SB = V_{m_k}^\dagger B = [B_p; 0_{(m_k-p)\times p}] \in \mathbb{R}^{m_k \times p}$. 

To find the maximizer of $\|S{\mathcal R_k}(s)\|$ efficiently, we perform a reduced QR factorization $Q_ZR_Z = S(AV_{m_k}-V_{m_k}T_{m_k}) \in \mathbb{R}^{s \times m_k}$, so that $$\|S{\mathcal R_k}(s)\| = \left\|R_Z(T_{m_k}-sI)^{-1}[B_p;\,0]\right\|.$$ Once the maximizer $s_k \in W$ of $\|S{\mathcal R_k}(s)\|$ has been found, we compute the leading $\ell_k$ right singular vectors of $S{\mathcal R_k}(s_k)$, referred to as $d_k \in \mathbb{R}^{p \times \ell_k}$. Then we determine a new pole $\widehat{s}_k$ to expand the rational Krylov subspace by the new vectors $U_k = (A-\widehat{s}_kI)^{-1}{\mathcal R_k}(s_k)d_k \in \mathbb{R}^{n \times \ell_k}$. Here, unlike the existing tangential RKSM taking $\widehat{s}_k = s_k$, the new pole $\widehat{s}_k$ may or may not be equal to the maximizer $s_k$. 
\begin{eqnarray}\nonumber
    U_k &=& (A-s_k I)^{-1}{\mathcal R_k}(s_k) = (A-\widehat{s}_kI)^{-1}(AV_{m_k}-V_{m_k}T_{m_k})(T_{m_k}-s_k I)^{-1}[B_p;\, 0]d_k \\ \nonumber
    &=& (A-\widehat{s}_k I)^{-1}\left[(A-\widehat{s}_k I)V_{m_k}-V_{m_k}(T_{m_k}-\widehat{s}_kI)\right](T_{m_k}-s_kI)^{-1}[B_p;\, 0]d_k \\ \nonumber
    &=& \left[V_{m_k}-(A-\widehat{s}_kI)^{-1}V_{m_k}(T_{m_k}-\widehat{s}_kI)\right](T_{m_k}-s_kI)^{-1}[B_p;\, 0]d_k.
\end{eqnarray}
Since $V_{m_k}(T_m-s_kI)^{-1}[B_p; 0_{(m_k-p)\times p}]d_k$ spans a subspace of $\mathrm{span}(V_{m_k})$, the new vectors $U_k$ can be taken as 
\begin{eqnarray}\label{ukdef1}
    U_k = (A-\widehat{s}_kI)^{-1} \left(V_{m_k}(T_{m_k}-\widehat{s}_kI)(T_{m_k}-s_kI)^{-1} [B_p;\, 0]d_k\right),
\end{eqnarray}
and if $\widehat{s}_k = s_k$, it is simplified as 
\begin{eqnarray}\label{ukdef2}
    U_k = (A-s_k I)^{-1} \left(V_{m_k}[B_p;\, 0]d_k\right) = (A-s_k I)^{-1}Bd_k. 
\end{eqnarray}

With the flexibility of choosing a new pole $\widehat{s}_k$ different from the maximizer $s_k$, it is possible to adopt a few repeated poles through sparse LU factorizations. At the first few steps, we may let $\widehat{s}_k = s_k$, perform and store the LU factorizations of $A-s_k I$; after these initial steps, each new pole $\widehat{s}_k$ can be chosen as the repeated pole closest to $s_k$ on a logarithmic scale. If an iterative linear solver is needed, this option may still make sense, since only a few different preconditioners will be constructed. The savings in runtime achieved by reusing the few sparse triangular factors or preconditioners might compensate the higher cost incurred by the slower convergence of RKSM with a few repeated poles than with adaptive poles ($\widehat{s}_k = s_k$). {We do not claim that sketching is the only possible route to repeated poles. A non-sketched TRKSM with a properly modified repeated-pole rational basis may also be possible, but it would require a separate reformulation of the basis construction, the residual and pole-selection machinery; this direction could be considered in a separate study.}

\subsection{Fast residual norm estimate}\label{subsect:resest} The residual norm estimate of the Lyapunov equation is another costly computation for the standard tangential RKSM; see section~\ref{sect:review}. With randomized sketching, this can also be done at a considerably lower cost. Let $T_{m_k}Y_{m_k}+Y_{m_k}T_{m_k}^T + [B_p;\,0][B_p;\,0]^T = 0$ be the projected Lyapunov equation right before step $k$, where $T_{m_k} = (SV_{m_k})^\dagger(SAV_{m_k})$, and the approximate solution is $X_{m_k} = V_{m_k}Y_{m_k}V_{m_k}^T$. The residual of the equation can be written as 
\begin{eqnarray}\label{lyapres}
    && AX_{m_k} + X_{m_k} A^T + BB^T = AX_{m_k} + X_{m_k} A^T + V_{m_k}[B_p;\,0][B_p;\,0]^TV_{m_k}^T \\ \nonumber
    &=& AV_{m_k}Y_{m_k}V_{m_k}^T + V_{m_k}Y_{m_k}V_{m_k}^TA^T - V_{m_k}(T_{m_k}Y_{m_k}+Y_{m_k}T_{m_k}^T)V_{m_k}^T \\ \nonumber
    &=& (AV_{m_k}-V_{m_k}T_{m_k})Y_{m_k}V_{m_k}^T + V_{m_k}Y_{m_k}(AV_{m_k}-V_{m_k}T_{m_k})^T. 
\end{eqnarray}
The norm of the residual can be estimated from the following result.

\begin{proposition}\label{prop:lyapres}
    Let $Q_ZR_Z = S(AV_{m_k}-V_{m_k}T_{m_k}) \in \mathbb{R}^{s \times m_k}$ and $Q_V R_V = SV_{m_k}\in \mathbb{R}^{s \times m_k}$ be reduced QR factorizations, and assume that the sketching matrix $S$ satisfies the sketching assumption \eqref{sketchassumption} for any $v \in \mathrm{span}(AV_{m_k}-V_{m_k}T_{m_k})$ and $v \in \mathrm{span}(V_{m_k})$. Then the following bound holds
    \begin{eqnarray}\label{lyapresnormest}
        &&\left\|AX_{m_k}+X_{m_k}A^T+BB^T\right\|_* \leq \frac{2}{1-\epsilon}\left\|R_ZY_{m_k}R_V^T\right\|_*,
    \end{eqnarray}
    where $\|\cdot\|_*$ can be the matrix 2-norm or the Frobenius norm.
\end{proposition}
\begin{proof}
    Given the two reduced QR factorizations, $Q_Z = S(AV_{m_k}-V_{m_k}T_{m_k})R_Z^{-1}$ and $Q_Z = SV_{m_k}R_V^{-1}$ are $s \times m_k$ matrices with orthonormal columns. From \eqref{lyapres}, 
    \begin{eqnarray}\nonumber
        &&\left\|AX_{m_k}+X_{m_k}A^T+BB^T\right\|_* \leq 2 \left\|(AV_{m_k}-V_{m_k}T_{m_k})Y_{m_k}V_{m_k}^T\right\|_* \\ \nonumber
        &=& 2\left\|(AV_{m_k}-V_{m_k}T_{m_k})R_Z^{-1}R_ZY_{m_k}R_V^T(V_{m_k}R_V^{-1})^T\right\|_* \\ \nonumber
        & \leq & 2 \left\|(AV_{m_k}-V_{m_k}T_{m_k})R_Z^{-1}\right\|_2\left\|R_ZY_{m_k}R_V^T\right\|_*\left\|V_{m_k}R_V^{-1}\right\|_2 \\ \nonumber
        & \leq & \frac{2}{1-\epsilon}\left\|S(AV_{m_k}-V_{m_k}T_{m_k})R_Z^{-1}\right\|_2\left\|SV_{m_k}R_V^{-1}\right\|_2 \left\|R_ZY_{m_k}R_V^T\right\|_* \\ \nonumber
        & = & \frac{2}{1-\epsilon} \|Q_Z\|_2\|Q_V\|_2\left\|R_ZY_{m_k}R_V^T\right\|_*= \frac{2}{1-\epsilon}\left\|R_ZY_{m_k}R_V^T\right\|_*,
    \end{eqnarray}
where the last inequality above is based on \eqref{F2normsineq} under the sketching assumption. Also, if $\|\cdot\|_*$ is the Frobenius norm, note that $\|FG\|_F \leq \|F\|_{p}\|G\|_F$ ($1 \leq p \leq \infty$) for any operator norm $\|\cdot\|_p \leq \|\cdot\|_F$. 
\end{proof}

Similar to the savings in the new pole selection by finding the norm maximizer of the sketched RKSM residual $S{\mathcal R_k}(s)$, the above proposition shows that the cost of estimating the residual norm of the Lyapunov equation can also be lowered significantly by computing the reduced QR factorizations of two $s \times m_k$ matrices, both of which are formed or stored explicitly. The bound in \eqref{lyapresnormest} is easy to compute since it only involves multiplications and norm evaluation of square matrices of order $m_k$. 

{We note that the residual norm bound \eqref{lyapresnormest} in Proposition~\ref{prop:lyapres} is a generalization of similar relations (3.5) and (3.6) in \cite{PalittaEtal25}, though our bound does not rely on any particular Arnoldi-like relations. In the special case of the sketched-and-truncated polynomial Arnoldi method, in our notation, the Arnoldi relation is $$AV_{m_k}=V_{m_k}T_{m_k}+\boldsymbol{v_{m_{k+1}}} \boldsymbol{h_{m_{k+1},m_k}}\boldsymbol{e_{m_k}^T},$$where $\boldsymbol{v_{m_{k+1}}} \in \mathbb{R}^{n \times \ell_k}$ with orthonormal columns with respect to the inner product defined by $S^TS$, $\boldsymbol{h_{m_{k+1},m_k}} \in \mathbb{R}^{\ell_k \times \ell_k}$, and $\boldsymbol{e_{m_k}^T} = [0\,\,I_{\ell_k}]\in \mathbb{R}^{\ell_k \times m_k}$. The relation \eqref{lyapres} leading to Proposition~\ref{prop:lyapres} shows that the two-sided sketched residual norm satisfies
\begin{eqnarray}\nonumber
&&\left\|S\left(AX_{m_k} + X_{m_k} A^T + BB^T \right)S^T\right\| \leq 2 \left\|S\left((AV_{m_k}-V_{m_k}T_{m_k})Y_{m_k}V_{m_k}^T\right)S^T\right\| \\ \nonumber
&=& 2\left\|S\left(\boldsymbol{v_{m_{k+1}}} \boldsymbol{h_{m_{k+1},m_k}}\boldsymbol{e_{m_k}^T} Y_{m_k}V_{m_k}^T\right)S^T\right\| = 2 \left\|\boldsymbol{h_{m_{k+1},m_k}}\boldsymbol{e_{m_k}^T} Y_{m_k}(SV_{m_k})^T\right\|.
\end{eqnarray}
If, in addition, $V_{m_k}$ contains orthonormal columns with respect to the inner product defined by $S^TS$, then the final expression above equals $2\left\|\boldsymbol{h_{m_{k+1},m_k}}\boldsymbol{e_{m_k}^T}Y_{m_k} \right\|$, which corresponds to (3.5) in \cite{PalittaEtal25} applied to the Lyapunov equation. The difference is that our bound in \eqref{lyapresnormest} is based on a one-sided sketched residual, whereas the bound in (3.6) of \cite{PalittaEtal25} is based on a two-sided sketched residual.
}

\subsection{Algorithm outline} With the details of new pole selection and residual estimation discussed, we provide a detailed description of the new tangential RKSM with randomized sketching. After the initialization (line 1), we expand and store $SAV_{m_k}$, $SV_{m_k}$, update the reduced QR factorization $Q_VR_V = SV_{m_k}$ and $T_{m_k} = (SV_{m_k})^\dagger SAV_{m_k}$, form $SAV_{m_k}-SV_{m_k}T_{m_k}$ and obtain its reduced QR factorization $Q_ZR_Z$ (line 3). Then we solve the projected Lyapunov equation (line 4), evaluate the equation residual norm (line 5) to see if the stopping criterion is satisfied (line 6-7). If the algorithm has not converged yet, the new ideal pole $s_k$ is determined as the maximizer of the sketched RKSM residual norm $\|S{\mathcal R_k}(s)\|$, and the tangential directions $d_k$ (leading singular vectors of $S{\mathcal R_k}(s_k)$) are computed (line 9). The new actual pole $\widehat{s}_k$ is then determined (depending on whether adaptive poles or a few repeated poles are adopted), and the new basis vectors $U_k$ are computed (line 10). Finally, $U_k$ is selectively orthogonalized against $V_{m_k}$ and the basis vectors are expanded to $V_{m_{k+1}}$, and the next step will be started (line 11). Failure of convergence is reported if the stopping criterion is not satisfied when the subspace dimension exceeds the predetermined maximum value (line 13-14). 

\begin{algorithm}
	\caption{A randomly sketched tangential RKSM for the Lyapunov equation}
	\label{alg:rstgrksm}
	\begin{algorithmic}[1]
	    \REQUIRE a stable matrix $A \in \mathbb{R}^{n \times n}$, $B \in \mathbb{R}^{n \times p}$ of rank $p$, relative tolerance $tol > 0$, maximum subspace dimension $m$, and sketching matrix $S \in \mathbb{C}^{s \times n}$
        \ENSURE $Z_{m_k} \in \mathbb{R}^{n \times m_k}$ ($m_k \leq m$) for the approximate solution $X_{m_k} = Z_{m_k} Z_{m_k}^T$;
        \STATE initialize a factorization $V_p B_p = B$ with $B_p \in \mathbb{R}^{p \times p}$; let $k=1$, $m_1 = p$.
        \WHILE {$\mathrm{dim}(V_{m_k}) < m$ and tolerance not satisfied}
        \STATE expand, update and store $SAV_{m_k}$, $SV_{m_k}$ and the reduced QR factorization $Q_VR_V = SV_{m_k}$; evaluate $T_{m_k}=(SV_{m_k})^\dagger SAV_{m_k}$, and obtain the reduced QR factorization $Q_ZR_Z = SAV_{m_k}-SV_{m_k}T_{m_k}$;
        \STATE solve the projected equation $T_{m_k}Y_{m_k}+Y_{m_k}T_{m_k}^T + [B_p;\,0][B_p;\,0]^T = 0$;
        \STATE evaluate the estimated residual norm $\eta_k = \frac{2}{1-\epsilon}\left\|R_ZY_{m_k}R_V^T\right\|$ in \eqref{lyapresnormest}; 
        \IF {the relative residual norm $\frac{\eta_k}{\|B^TB\|} \leq tol$}
            \STATE define $X_{m_k}= Z_{m_k}Z_{m_k}^T$, where $Z_{m_k} = V_{m_k}L_{m_k}$ and $L_{m_k}L_{m_k}^T = Y_{m_k}$ is the Cholesky factorization of $Y_{m_k}$; return $Z_{m_k}$;
        \ENDIF
        \STATE evaluate $s_k = \mathrm{argmax}_{s \in W}\|S{\mathcal R_k}(s)\| = \left\|R_Z(T_{m_k}-sI)^{-1}[B_p;\,0_{(m_k-p)\times p}]\right\|$, and find the leading $\ell_k$ right singular vectors of $S{\mathcal R_k}(s_k)$, assembled in $d_k \in \mathbb{R}^{p \times \ell_k}$; 
        \STATE determine a new pole $\widehat{s}_k$, and evaluate $U_k \in \mathbb{R}^{n \times \ell_k}$ by \eqref{ukdef1} (if $\widehat{s}_k \ne s_k$, e.g., using a few repeated poles) or \eqref{ukdef2} (if $\widehat{s}_k = s_k$, e.g., using adaptive poles);
        \STATE let $m_{k+1} = m_k + \ell_k$, and apply selective orthogonalization of $U_k$ against $V_{m_k}$ and append the updated vectors to $V_{m_k}$ to obtain $V_{m_{k+1}}$; let $k = k+1$;
        \ENDWHILE
        \IF {tolerance not achieved}
        \STATE declare failure of convergence and return the most recent $Z_{m_k} = V_{m_k}L_{m_k}$;
        \ENDIF
    \end{algorithmic} 
\end{algorithm}

\section{A computable error bound of approximate solutions} In this section, we first establish a property of exactness of the sketched RKSM, which is well known for the standard RKSM. We then provide a computable bound showing the dependence of the error of the approximate solution on the pole selection of RKSM. {The key idea is to (i) construct an infinite series of rational functions, depending on a parameter $t \geq 0$, to represent the action of matrix exponentials $e^{tA}b$, and substitute it into the exact solution $X= \int_0^\infty e^{tA}b\left(e^{tA}b\right)^*dt$ of the Lyapunov equation, assuming that $A$ is stable; (ii) construct a similar infinite series for the approximate solution $X_m$ obtained from the $(m-1)$-step RKSM; (iii) use the property of exactness (see Theorem~\ref{thm:exactrksm}) to show that the two infinite series share the common leading $(m-1)^2$ terms, so that the error $\|X-X_m\|$ can be bounded by the remainder.}

\subsection{Property of exactness} {The standard polynomial and rational Arnoldi methods are widely adopted for approximating functions of matrices times vectors $f(A)b$ or solving the Lyapunov equation $AX+XA^T+bb^T=0$. Some error analyses in both settings depend on the property of exactness; see, e.g., \cite[Theorem~4.3]{Cortinovis.etal.2024} and \cite[Corollary~2.3]{GS23} for the polynomial case, as well as \cite[Lemma~3.1]{Gue13b} and \cite[Lemma~2.2]{Xu.Xue.2023} for the rational case.} For rational Arnoldi, let $V_m \in \mathbb{R}^{n \times m}$ contain {an} orthonormal basis of the rational Krylov subspace $\mathcal{Q}_m(A,b) = q_{m-1}(A)^{-1}\mathcal{K}_m(A,b)$, where $q_{m-1}(z) = \prod_{k=1}^{m-1}(z-s_k)$ is a polynomial of degree $m-1$ whose roots are the poles used to generate $\mathcal{Q}_m(A,b)$, and $\mathcal{K}_m(A,b)=\mathrm{span}\{b,Ab,\ldots,A^{m-1}b\}$ is the standard Krylov subspace of dimension $m$. The property of exactness states that $$r_{m-1}(A)b = V_m r_{m-1}(H_m)e_1\|b\|,$$where $H_m = V_m^T A V_m$, and $r_{m-1}(z)$ is any rational function of $(j,m-1)$-type ($0 \leq j \leq m-1$) with the poles $s_1,\ldots,s_{m-1}$ used to generate $\mathcal{Q}_m(A,b)$. 

We now show that this relation also holds for the sketched RKSM with non-orthonormal basis $V_m$ and the sketched compressed matrix $T_m = (SV_m)^\dagger SAV_m$. The full statement is given in the following theorem. 

\begin{theorem}\label{thm:exactrksm}
    Given $A \in \mathbb{R}^{n \times n}$, let $b \in \mathbb{R}^{n}$ be the initial vector, and $s_1,\ldots,s_{m-1}$ the poles to generate the rational Krylov subspace $\mathcal{Q}_m(A,b) = q_{m-1}(A)^{-1}\mathcal{K}_m(A,b)$, where $q_{m-1}(z) = \prod_{k=1}^{m-1}(z-s_k)$. Let $V_m\in \mathbb{R}^{n \times m}$ contain a general basis of $\mathcal{Q}_m(A,b)$ with $V_me_1 = v_1 = b$ in direction, $S \in \mathbb{C}^{s \times n}$ be a full-rank sketching matrix, and assume that $SV_m$ is of full rank $m$. Define the sketched compressed matrix $T_m = (SV_m)^\dagger SAV_m$, and let $r_{m-1}(z)$ be any rational function of $(j,m-1)$-type, where $0 \leq j \leq m-1$ and the poles of $r_{m-1}(z)$ are those used to generate $\mathcal{Q}_m(A,b)$. Then
    \begin{eqnarray}\label{exactprop}
        r_{m-1}(A)b = V_mr_{m-1}(T_m)e_1 \|b\|/\|v_1\|. 
    \end{eqnarray}
\end{theorem}
\begin{proof}
    Define $f=q_{m-1}(A)^{-1}b$, such that any vector $g \in \mathcal{Q}_m(A,b)$ can be written as $g = p_{m-1}(A)q_{m-1}(A)^{-1}b=p_{m-1}(A)f$. Therefore, for any index $0 \leq j \leq m-1$, $A^j f \in \mathrm{span}(V_m)$, so that there exists $u\in \mathbb{R}^m$ satisfying $V_m u = A^j f$, and hence $V_m V_m^\dagger A^j f = V_m u = A^j f$. Meanwhile, we note that $SA^jf \in \mathrm{span}(SV_m)$, so that $u = (SV_m)^\dagger (SA^j f) = V_m^\dagger A^j f$, and it follows that $V_m (SV_m)^\dagger (SA^j f) = V_m u =  A^j f$. 

    Note that $V_m V_m^\dagger f = f$ because $f \in \mathrm{span}(V_m)$. It follows that
    \begin{eqnarray}\label{sketchrk}
    && A^j f = V_m u = V_m (SV_m)^\dagger (SA^jf) = V_m (SV_m)^\dagger SA \left( A^{j-1}f\right) \\ \nonumber
    &=& \left(V_m (SV_m)^\dagger SA\right) \left( V_m (SV_m)^\dagger SA \right) A^{j-2}f = \ldots = \left(V_m (SV_m)^\dagger SA \right)^j V_m V_m^\dagger f\\ \nonumber
    &=& V_m \left((SV_m)^\dagger (SAV_m)\right)^j (V_m^\dagger f) = V_m T_m^j (V_m^\dagger f).
    \end{eqnarray}

    Since $q_{m-1}(A)$ is a linear combination of $A^j$'s of different powers $j$ ($0 \leq j \leq m-1$), \eqref{sketchrk} leads to $q_{m-1}(A)f = V_m q_{m-1}(T_m)\left(V_m^\dagger f\right)$. Multiplying both sides of this relation on the left by $V_m^\dagger$ and then by $q_{m-1}(T_m)^{-1}$, we have
    \begin{eqnarray}\label{vmpinvf}
        V_m^\dagger f = q_{m-1}(T_m)^{-1}V_m^\dagger q_{m-1}(A)f = q_{m-1}(T_m)^{-1}V_m^\dagger b.
    \end{eqnarray}
    
    It then follows from \eqref{sketchrk} and \eqref{vmpinvf} that 
    \begin{eqnarray}
        && A^j q_{m-1}(A)^{-1}b = A^j f =  V_m T_m^j (V_m^\dagger f) \\ \nonumber
        &=& V_m T_m^j q_{m-1}(T_m)^{-1}V_m^\dagger b = V_m T_m^j q_{m-1}(T_m)^{-1}V_m^\dagger b.
\end{eqnarray}
    Since $p_{m-1}(A)$ is a linear combination of $A^j$ of different power $j$, it follows that 
    \begin{eqnarray}
        r_{m-1}(A)b &=& p_{m-1}(A)q_{m-1}(A)^{-1}b = V_m p_{m-1}(T_m)q_{m-1}(T_m)^{-1}V_m^\dagger b \\ \nonumber
        &=& V_mr_{m-1}(T_m)e_1 \|b\|/\|v_1\|,
    \end{eqnarray}
    where the last equality follows from the fact that $b = v_1$ in direction. 
\end{proof}

\subsection{Solutions as in infinite series of rational functions} We show through the Faber-Dzhrbashyan (FD) rational function and the property of exactness (Theorem \ref{thm:exactrksm}) that the convergence of the randomly sketched RKSM is essentially the same as that of the standard RKSM shown in \cite{Druskin.etal.2011}. We provide a bound showing an explicit dependence of the error of the approximate solution on the poles of RKSM. 

\subsubsection{FD series form of $f(A)b$} We first review the Faber-Dzhrbashyan (FD) rational functions for our error estimate. Let $E\subset \mathbb{C}$ be a connected compact domain in the complex plane, $D = \{w\,:\,|w|\leq 1\}$ be the unit disk, $\phi: \overline{\mathbb{C}}\setminus{E} \rightarrow \overline{\mathbb{C}}\setminus{D}$ be the Riemann mapping and $\psi = \phi^{-1}$ be the inverse mapping. Let $\left\{z_k\right\}_{k=1}^{\infty}$ be a sequence of points lying in the exterior of $E$, and let $\alpha_k=1/\overline{\phi(z_k)}\in D = \{w:|w| \le 1\}$ ($k\ge 1$). The Takenaka-Malmquist (TM) system of rational functions are defined as
\begin{eqnarray} \label{eq:varphi}
  \varphi_1(w) = \frac{\sqrt{1-|\alpha_1|^2}}{1-\overline{\alpha_1}w}, \quad
  \varphi_j(w) = \frac{\sqrt{1-|\alpha_j|^2}}{1-\overline{\alpha_j}w}\prod_{k=1}^{j-1}\frac{\alpha_k-w}{1-\overline{\alpha_k}w} \frac{|\alpha_k|}{\alpha_k}\,\ (j\ge 2),
\end{eqnarray}
where ${|\alpha_k|}/{\alpha_k}=1$ if $\alpha_k = 0$. Since $\phi(z_k) \in \overline{\mathbb{C}}\setminus{D}$, $|\alpha_k| = \big|1/\overline{\phi(z_k)}\big|< 1$. The TM systems $\left\{\varphi_n(w)\right\}_{n=0}^{\infty}$ form an orthonormal basis on the unit circle, i.e., $\langle \varphi_m(w),\varphi_n(w) \rangle := \frac{1}{2\pi}\int_{0}^{2\pi}\varphi_m(e^{it})\overline{\varphi_n(e^{it})}dt = \delta_{mn}$, where $\delta_{mn}$ is the Kronecker delta. 

The FD rational function $M_j(z)$ ($j \ge 1$) is defined as the summation of the principal part in the Laurent series of $\varphi_j(\phi(z))$ centered at the points $\left\{z_k\right\}_{k=1}^{j}$. Therefore, it can be written as $M_j(z)=\frac{p_{j-1}(z)}{\prod_{k=1}^{j}(z-z_k)}$, where $p_{j-1}(z)$ is a polynomial with degree $\leq j-1$. Let $\Gamma_D$ be the preimage of the unit disk under the map $w=\phi(z)$ and $G_D$ be the interior of $\Gamma_D$, the FD rational functions can also be written in the form (see \cite[Section 13, equation (4)]{Suetin.1998}):
\begin{equation}\label{eq:ddf2_mjz}
  M_j(z)=\frac{1}{2\pi i}\int_{\Gamma_D}\frac{\varphi_{j}\left(\phi(\zeta)\right)}{\zeta-z}d\zeta,~z\in G_D.
\end{equation}

If the conditions $\sum_{k=1}^{\infty}\left(1-|\alpha_k|\right)=\infty$ and $\lim_{r\to 1^+}\int_{0}^{2\pi}\left|\psi'(re^{i\theta})\right|^2d\theta <\infty$ hold, we have the expansion $\frac{\psi'(w)}{\psi(w)-z} = \frac{1}{w}\sum_{j=1}^{\infty}\overline{\varphi_{j}\left(\frac{1}{\overline{w}}\right)}M_j(z)$ with $z\in G$ and $|w|>1$, where $G$ is the interior of $E \supset W(A)$; see, e.g., \cite[page 259]{Suetin.1998}. If $W(A)\subset G$, then
\begin{equation}\label{eq:psipar}
  \psi'(w)\left(\psi(w)I-A\right)^{-1} = \frac{1}{w}\sum_{j=1}^{\infty}\overline{\varphi_{j}\left(\frac{1}{\overline{w}}\right)}M_j(A),~ |w|>1, 
\end{equation}
where upper bounds on $\left|\varphi_j(1/\overline{w})\right|$ and $\|M_j(A)\|$ are given in the following lemma.

\begin{lemma}\label{lemma:maxvarphi}
Given $\tau \geq 1$, and $|\alpha_j|<1$ for all $j\ge 1$, it holds that
\begin{eqnarray}\label{eq:lemmavarphi}
 && \|M_j(A)\| \leq 2\sqrt{\frac{1+|\alpha_{j}|}{1-|\alpha_{j}|}} + \sqrt{1-|\alpha_{j}|^2}\prod_{k=1}^{j-1}\left|\alpha_k\right|, \quad\mbox{and}  \\ \nonumber
 && \max_{|w|=\tau}\left|\overline{\varphi_j\left(\frac{1}{\overline{w}}\right)}\right| \leq \frac{\tau\sqrt{1-|\alpha_j|^2}}{\tau-|\alpha_j|} \prod_{k=1}^{j-1}\frac{\tau|\alpha_k|+1}{\tau+|\alpha_k|}. 
\end{eqnarray}
\end{lemma}
\begin{proof}
    See \cite[Appendix]{Xu.Xue.2023} for details. 
\end{proof}

Define $E_{\tau}=E\cup\left\{z~\arrowvert~z\in \mathbb{C}\setminus E, |\phi(z)|\le \tau\right\}$ ($\tau \geq 1$) and let $\Gamma_{\tau}$ be the boundary of $E_\tau$. If $f(z)$ is analytic in $E_{\tau}$, since $W(A)\subset E \subset E_{\tau}$, by (\ref{eq:psipar}), it holds that
\begin{eqnarray}\nonumber
  f(A) &=& \frac{1}{2\pi i}\int_{\Gamma_{\tau}}f(z)(zI-A)^{-1}dz
  = \frac{1}{2\pi i}\int_{|w|=\tau}f\left(\psi(w)\right)\left(\psi(w)I-A\right)^{-1}\psi'(w)dw \\ \label{eq:faexpan}
   &=& \sum_{j=1}^{\infty}M_j(A) \frac{1}{2\pi i}\int_{|w|=\tau}\frac{f\left(\psi(w)\right)}{w} \overline{\varphi_{j}\left(\frac{1}{\overline{w}}\right)} dw = \sum_{j=1}^{\infty}c_jM_j(A),
\end{eqnarray}
where the coefficient 
\begin{eqnarray}\label{fdcoeffs}
c_j &=& \frac{1}{2\pi i}\int_{|w|=\tau}\frac{f\left(\psi(w)\right)}{w} \overline{\varphi_{j}\left(\frac{1}{\overline{w}}\right)} dw = \frac{1}{2\pi}\int_0^{2\pi} f\left(\psi(\tau e^{i\theta})\right)\overline{\varphi_j(e^{i\theta}/\tau)}d\theta  \\ \nonumber
&=& \frac{1}{2\pi i}\int_{|w|=\tau} \frac{f(\psi(w))}{w}\frac{w\sqrt{1-|\alpha_j|^2}}{w-\alpha_j}\prod_{i=1}^{j-1}\frac{w\overline{\alpha_i}-1}{w-\alpha_i}\frac{|\alpha_i|}{\overline{\alpha_i}}dw.
\end{eqnarray}
Since $w=0,\alpha_1,\ldots,\alpha_j$ are the poles of the integrand, the value of $c_j$ is \emph{independent} of $\tau \geq 1$ as long as $f$ is analytic in $E_\tau$. 

Equations \eqref{eq:faexpan} and \eqref{fdcoeffs} define the expansion of $f(A)$ in the FD rational basis functions. For a large matrix $A \in \mathbb{R}^{n \times n}$, the goal of approximating $f(A)b$ numerically by RKSM in $m$ steps ($m \ll n$) is to find $m$ nearly optimal poles $\{z_k\}_{k=1}^{m}$ (which define $\{\alpha_k\}_{k=1}^{m}$) such that the partial sum $\sum_{j=1}^{m}c_jM_j(A)b \approx f(A)b$ is accurate. 

\subsubsection{Error as the remainder of FD series} In this section, we derive a computable upper bound on the error norm of the approximate solution generated by the standard or the sketched RKSM. By \eqref{eq:faexpan} and \eqref{fdcoeffs}, the exact solution of the Lyapunov equation $AX+XA^T+bb^T= 0$ for a stable matrix $A$ is
\begin{eqnarray}\nonumber
    X &=& \int_0^{\infty}e^{tA}b \big(e^{tA}b\big)^* dt = \int_0^{\infty} \bigg(\sum_{k=1}^\infty c_k(t) M_k(A)b\bigg)\bigg(\sum_{\ell=1}^\infty c_\ell(t)M_\ell(A)b\bigg)^* dt \\ \label{exactsol}
    &=& \sum_{k=1}^\infty \sum_{\ell=1}^\infty \left(\int_0^{\infty}c_k(t)\overline{c_\ell(t)}dt\right) M_k(A)b \left(M_\ell(A)b\right)^*,  \\ \nonumber
    && \mbox{ where } c_k(t) =  \frac{1}{2\pi}\int_0^{2\pi}e^{t\psi\left(\tau e^{i\theta}\right)}\overline{\varphi_{k}\left({e^{i\theta}}/{\tau}\right)}d\theta.
\end{eqnarray}
The coefficients with double running indices in \eqref{exactsol} can be written as
\begin{eqnarray}\label{betakl}
    && \beta_{k,\ell} = \int_0^{\infty} c_k(t)\overline{c_\ell(t)} dt \\ \nonumber
&=& \int_0^{\infty}  \left(\frac{1}{2\pi}\int_0^{2\pi}e^{t\psi\left(\tau e^{i\theta}\right)}\overline{\varphi_{k}\left({e^{i\theta}}/{\tau}\right)}d\theta\right) \left( \frac{1}{2\pi}\int_0^{2\pi}e^{t\overline{\psi\left(\tau e^{i\eta}\right)}}{\varphi_{\ell}\left({e^{i\eta}}/{\tau}\right)}d\eta\right)  dt \\ \nonumber
&=& \frac{1}{4\pi^2}\int_0^{2\pi}\int_0^{2\pi} \int_0^{\infty}e^{t\left(\psi\left(\tau e^{i\theta}\right)+\overline{\psi\left(\tau e^{i\eta}\right)}\right)}dt \, \overline{\varphi_{k}\left({e^{i\theta}}/{\tau}\right)}\, {\varphi_{\ell}\left({e^{i\eta}}/{\tau}\right)}d\eta d\theta \\ \nonumber
&=& -\frac{1}{4\pi^2}\int_0^{2\pi}\int_0^{2\pi} \frac{1}{\psi\left(\tau e^{i\theta}\right)+\overline{\psi\left(\tau e^{i\eta}\right)}}\overline{\varphi_{k}\left({e^{i\theta}}/{\tau}\right)}\,{\varphi_{\ell}\left({e^{i\eta}}/{\tau}\right)}d\eta d\theta.
\end{eqnarray}
Here, $\tau \geq 1$ is not too large, so that $\max_{\theta,\eta \in [0,2\pi)}\mathrm{Re}\left(\psi(\tau e^{i\theta})+\overline{\psi(\tau e^{i\eta})}\right) < 0$ and the inner integral of the variable $t$ can be simplified. 

Meanwhile, the approximate solution obtained from an $(m-1)$-step standard or sketched RKSM is $X_m = V_m Y_m V_m^T$, where $V_m$ contains the basis vectors of the rational Krylov subspace {$\mathcal{Q}_m(A,b)$, and $Y_m$ is the solution to the projected (standard or sketched) Lyapunov equation $T_m Y + Y T_m^T + \frac{\|b\|^2}{\|v_1\|^2}e_1 e_1^T = 0$. We note that for the non-sketched RKSM with orthonormal basis $V_m$, the matrix $T_m = (SV_m)^\dagger SAV_m$ is replaced with $H_m = V_m^TAV_m$, and the results below still hold.} From \eqref{exactsol}, we have 
\begin{eqnarray}\label{approxsoln}
    && X_m = V_m Y_m V_m^T = V_m \left(\int_0^\infty \left(e^{tT_m}e_1\frac{\|b\|}{\|v_1\|}\right) \left(e^{tT_m}e_1\frac{\|b\|}{\|v_1\|}\right)^* dt \right) V_m^T \\ \nonumber
    &=& \sum_{k=1}^\infty \sum_{\ell=1}^\infty \int_0^\infty c_k(t)\overline{c_\ell(t)}dt \left(V_mM_k(T_m)e_1\frac{\|b\|}{\|v_1\|}\right)\left(V_mM_\ell(T_m)e_1\frac{\|b\|}{\|v_1\|}\right)^*.
\end{eqnarray}
Here, note that each $M_k(z)$ ($1 \leq k \leq m-1$) is a rational function of $(k-1,k)$-type, and therefore also is of $(m-1,m-1)$-type. As a result, it holds that $$M_k(A)b = V_m M_k(T_m)e_1 \frac{\|b\|}{\|v_1\|}\quad (1 \leq k \leq m-1)$$ by \eqref{exactprop} in Theorem \ref{thm:exactrksm}, so that the two sets of $(m-1)^2$ terms with $1\leq k \leq m-1$ and $1 \leq \ell \leq m-1$ in $X_m-X$ are canceled. Without loss of generality, assume that $\|b\|=1$ and let $\gamma_m = \frac{\|V_m\|^2}{\|v_1\|^2}$. It follows that the error $\|X_m-X\|$ can be written as
{\small
\begin{eqnarray}\label{errboundabstract}
&& \sum_{\max\{k,\ell\} \geq m} \left\|\beta_{k,\ell} \left(M_k(A)b\right)\left(M_\ell(A)b\right)^* - \beta_{k,\ell}\frac{\|b\|^2}{\|v_1\|^2}\left(V_mM_k(T_m)e_1\right)\left(V_mM_\ell(T_m)e_1\right)^*\right\| \\ \nonumber
& \leq & \left\{\sum_{k=m}^\infty \sum_{\ell = 1}^{m-1} +  \sum_{k=1}^{m-1} \sum_{\ell = m}^\infty + \sum_{k=m}^\infty \sum_{\ell = m}^\infty \right\}\big|\beta_{k,\ell}\big| \Big(\|M_k(A)\|\|M_\ell(A)\|+\gamma_m\|M_k(T_m)\|\|M_\ell(T_m)\|\Big).
\end{eqnarray}
}

For a given step $m-1$, we can set the poles that have never been used as $z_k = \infty$ ($k \geq m$), so that $\alpha_k = 1/\overline{\phi(z_k)} = 0$. It follows from \eqref{eq:lemmavarphi} that $\|M_k(A)\| \leq 2+1 = 3$. Due to the symmetry of $k$ and $\ell$, we have
{\small
\begin{eqnarray}\label{errbound}
    && \qquad \|X_m-X\| \leq \sum_{k=1}^{m-1} \sum_{\ell = m}^\infty 6(1+\gamma_m)\big|\beta_{k,\ell}\big|\|M_k(A)\| + \sum_{k=m}^\infty \sum_{\ell = m}^\infty 9\big|\beta_{k,\ell}\big|(1+\gamma_m) \\ \nonumber
    &\leq& (1+\gamma_m)\left(\sum_{k=1}^{m-1} \sum_{\ell = m}^\infty 6\big|\beta_{k,\ell}\big| \bigg(2\sqrt{\frac{1+|\alpha_k|}{1-|\alpha_k|}}+\sqrt{1-|\alpha_k|^2}\prod_{j=1}^{k-1}|\alpha_j| \bigg) + \sum_{k=m}^\infty \sum_{\ell = m}^\infty 9\big|\beta_{k,\ell}\big|\right). 
\end{eqnarray}
}The expression of the coefficients $\beta_{k,\ell}$ ($\ell \geq m$) in \eqref{errbound} depends on $k$. For $1\leq k \leq m-1$, \eqref{betakl} is written as
\begin{eqnarray}\nonumber
\beta_{k,\ell} &=& -\frac{1}{4\pi^2}\int_0^{2\pi} \int_0^{2\pi} \frac{(-1/\tau)^{\ell-m}}{\psi(\tau e^{i\theta})+\overline{\psi(\tau e^{i\eta})}}\prod_{j=1}^{m-1}\frac{\tau \alpha_j - e^{-i\eta}}{\tau - \overline{\alpha_j}e^{-i\eta}}\frac{|\alpha_j|}{\alpha_j} e^{-i(\ell-m)\eta}d \eta \times \\ \label{betaklA}
&& \qquad \frac{\tau\sqrt{1-|\alpha_k|^2}}{\tau-\alpha_k e^{-i\theta}} \prod_{j=1}^{k-1} \frac{\tau \overline{\alpha_j}-e^{-i\theta}}{\tau - \alpha_j e^{-i\theta}}\frac{|\alpha_j|}{\overline{\alpha_j}}d\theta \\ \nonumber
&=& \frac{(-1)^{\ell-m+1}}{4\pi^2}\int_0^{2\pi} \,\int_0^{2\pi} \frac{e^{-i(\ell-1)\eta}}{\psi(e^{i\theta})+{\psi(e^{-i\eta})}}\prod_{j=1}^{m-1}\left(\frac{e^{i\eta}\alpha_j - 1}{1 - \overline{\alpha_j}e^{-i\eta}}\frac{|\alpha_j|}{\alpha_j}\right)d \eta\, \times \\ \label{betaklAp}
&& \qquad e^{-i(k-1)\theta}\frac{\sqrt{1-|\alpha_k|^2}}{1-\alpha_k e^{-i\theta}} \prod_{j=1}^{k-1} \left(\frac{e^{i\theta}\overline{\alpha_j}-1}{1 - \alpha_j e^{-i\theta}}\frac{|\alpha_j|}{\overline{\alpha_j}}\right)d\theta. \quad (1\leq k \leq m-1, \ell \geq m)
\end{eqnarray}
In the last step above, we let $\tau = 1$ and split off a factor of $e^{-i\theta}$ or $e^{-i\eta}$ from each factor of the two Blaschke products. This highlights that both the numerators $e^{i\eta}-1$ (or $e^{i\theta}\overline{\alpha_j}-1$) and the denominator $1-\overline{\alpha_j}e^{-i\eta}$ (or $1-\alpha_j e^{-i\theta}$) in \eqref{betaklAp} are uniformly bounded away from the origin, since all $|\alpha_j|$'s are uniformly bounded below and away from $1$ by assumption. Therefore, the integrands of the inner integral of $\eta$ and the outer integral of $\theta$ are both analytic periodic functions. In addition, we assume that the inverse Riemann mapping satisfies $\psi(\overline{z}) = \overline{\psi(z)}$, so that $\overline{\psi(e^{i\eta})} = \psi(e^{-i\eta})$. 

Similarly, for $k \geq m$, \eqref{betakl} is written as
\begin{eqnarray}\label{betaklB}
\qquad    \beta_{k,\ell} &=& -\frac{1}{4\pi^2}\int_0^{2\pi}\int_0^{2\pi}\frac{(-1/\tau)^{k+\ell-2m}}{\psi(\tau e^{i\theta})+\overline{\psi(\tau e^{i\eta})}}e^{-i(\ell-m)\eta}e^{-i(k-m)\theta} \times \\ \nonumber
    && \qquad\qquad \prod_{j=1}^{m-1}\frac{\tau \alpha_j - e^{-i\eta}}{\tau - \overline{\alpha_j}e^{-i\eta}} \prod_{j=1}^{m-1}\frac{\tau \overline{\alpha_j} - e^{-i\theta}}{\tau - {\alpha_j}e^{-i\theta}}d\eta d\theta, \quad (k\geq m, \ell \geq m) \\ \nonumber
   &=& \frac{(-1)^{k+\ell-2m+1}}{4\pi^2}\int_0^{2\pi}\int_0^{2\pi}\frac{1}{\psi(e^{i\theta})+{\psi(e^{-i\eta})}}e^{-i(\ell-1)\eta}e^{-i(k-1)\theta} \times \\ \nonumber
    && \qquad\qquad \prod_{j=1}^{m-1}\frac{e^{i\eta}\alpha_j - 1}{1 - \overline{\alpha_j}e^{-i\eta}} \prod_{j=1}^{m-1}\frac{e^{i\theta}\overline{\alpha_j} - 1}{1 - {\alpha_j}e^{-i\theta}}d\eta d\theta
\end{eqnarray}

Since the integrands in the double/repeated integrals \eqref{betaklAp} and \eqref{betaklB} are analytic periodic functions, it is well known by the decay rate of Fourier series that there is a $\rho \in (0,1)$, depending on $\{\alpha_j\}_{j=1}^{m-1}$,  such that $|\beta_{k,\ell}| \leq \mathcal{O}(\rho^{k+\ell})$; see, e.g., \cite[Ch.\,8]{Korner.book.1988} and \cite[Ch.\,3]{Grafakos.book.2008}. In addition, since $|\alpha_k|$'s are bounded away from $1$, it follows from \eqref{errbound} and the infinite series $\sum_{q=1}^\infty q\rho^{m+q} = \frac{\rho^{m+1}}{(1-\rho)^2}$ that $\|X_m-X\| \leq \mathcal{O}(\rho^m)$. The above analysis is summarized in the following theorem. { To see a convergence analysis for sketched polynomial Arnoldi for Lyapunov equations based on Faber polynomials, we refer the readers to \cite[Theorem 4.6 \& Corollary 4.7]{PalittaEtal25}.}

{
\begin{theorem}
    Given a stable matrix $A \in \mathbb{R}^{n \times n}$, let $E\subset \mathbb{C}$ be a domain lying strictly in the left half complex plane, such that $W(A) \subset E$ and the Riemann mappings $\phi$ and $\psi = \phi^{-1}$ can be defined. For the Lyapunov equation $AX+XA^T+bb^T=0$, let $X$ and $X_m$ be the exact solution and the approximate solution from the $(m-1)$-step standard or sketched RKSM, with $\{z_k\}_{k=1}^{m-1} \in \overline{\mathbb{C}}\setminus E$ being the poles used by RKSM. For the sketched RKSM, assume that $T_m=(SV_m)^\dagger SAV_m$ satisfies $W(T_m) \subset E$. Then the error $\|X_m-X\|$ is bounded as in \eqref{errbound}, with $\beta_{k,\ell}$'s defined in \eqref{betaklA} and \eqref{betaklB}, where $\alpha_k = 1/\overline{\phi(z_k)}$, $\gamma_m = \|V_m\|/\|v_1\|$, and $\tau \geq 1$ satisfies $\max_{\theta,\eta \in [0,2\pi)}\mathrm{Re}\big(\psi(\tau e^{i\theta})+\overline{\psi(\tau e^{i\eta})}\big) < 0$. Suppose that $|\alpha_j|$'s \textup{($1 \leq j \leq m-1$)} are uniformly bounded away from $1$. Then there exists a $\rho \in (0,1)$ that depends on $\{\alpha_j\}_{j=1}^{m-1}$, such that $\|X_m-X\| \leq \mathcal{O}(\rho^m)$.
\end{theorem}
}

{\bf Discussions.} Equations \eqref{errbound}, \eqref{betaklA} and \eqref{betaklB} provide a bound on $\|X_m-X\|$ that depends on the mapped poles $\{\alpha_k\}_{k=1}^{m-1}$ ($\alpha_k = 1/\overline{\phi(z_k)}$). This bound is derived on the basis of the property of exactness, which holds for the standard and the sketched variant. Compared to the estimates in \cite{Knizhnerman.Simoncini.2011} and \cite{Druskin.etal.2011}, our expressions of the coefficients $\beta_{k,\ell}$'s and the bound on $\|M_j(A)\|$ are more specific, without the use of the $\lesssim$ or $\mathcal{O}$ notation. In particular, the bound on $\|M_j(A)\|$ \eqref{eq:lemmavarphi}  shows that the poles need to be uniformly separated from the boundary of $E$ to secure a uniform bound, and \eqref{betaklA} and \eqref{betaklB} give the definitions of $\beta_{k,\ell}$ on the unit circle in the $w$-plane, which do not rely on the derivative function $\psi'(\phi(z))$ and the contour homotopic in $\overline{\mathbb{C}} \setminus [E \cup (-E )]$ to the imaginary axis in the $z$-plane, since we integrated the scalar-valued exponential function in \eqref{betakl} instead of converting it into integrals along the imaginary axis.

More importantly, the coefficients $\beta_{k,\ell}$'s and the explicit bound \eqref{errbound} can be computed numerically, and their validity can be verified by comparing them with established results. We note that it would be quite costly to evaluate each $\beta_{k,\ell}$ separately by a quadrature for double integrals, as we may need tens of thousands of such coefficients for $1\leq k \leq m-1$ and millions of them for $k \geq m$. Fortunately, the form of the expressions in \eqref{betaklA} and \eqref{betaklB} strongly suggests that the integrals of periodic functions can be approximated by trapezoidal rules and evaluated by the fast Fourier transform (FFT). In fact, the inner integral of $\eta$ in \eqref{betaklA} can be approximated by a one-dimensional FFT for $m \leq \ell \leq m+N-1$, whereas the double integral in \eqref{betaklB} can be estimated by a two-dimensional FFT directly for $m \leq k,\ell \leq m+N-1$, where $2^p N$ is the size of input arrays (number of uniform subintervals of $[0,2\pi)$) for FFTs, and we may need $2^p \geq 8$ to avoid undersampling and inaccurate quadrature. 

{
For the error bound \eqref{errbound} to reflect the observed convergence of RKSM, we take $E$ to be the smallest domain strictly contained in the open left half-plane that contains both $W(A)$ and $W(T_m)$. This requirement ensures that \eqref{exactsol}--\eqref{approxsoln} are applicable, and hence that \eqref{errboundabstract} and \eqref{errbound} are valid. In practice, for a given $m$ such a domain may fail to exist, since $W(T_m)$ (with $T_m=(SV_m)^\dagger SAV_m$ rather than $H_m=V_m^TAV_m$) need not be contained in $W(A)$ and may cross the imaginary axis. Accordingly, we assume that $W(T_m)\subset E$ so that the bound applies; if this condition is violated, \eqref{errbound} does not provide a guarantee for $\|X_m-X\|$.
}{Nevertheless, this is only a \emph{sufficient} condition for the FD-based bound, not a solvability condition for the projected Lyapunov equation. The latter still has a unique algebraic solution whenever \(\lambda_i(T_m)+\lambda_j(T_m)\ne0\) for all \(i,j\). If \(W(T_m)\) crosses the imaginary axis, the particular bound \eqref{errbound} is not certified, but this does not imply failure of the sketched method or exclude other possible error estimates; see, e.g., \cite[Section 5.2.1]{Simoncini2016}.}

{
Choosing $E$ as small as possible is essential for sharpness. If one replaces the minimal admissible $E$ by a larger set still containing $W(A)$ and $W(T_m)$ while keeping the RKSM poles $\{z_k\}$ fixed outside $E$, then the poles become closer to the boundary of the enlarged domain. Consequently, the mapped poles $\alpha_k=1/\overline{\phi(z_k)}$ approach the unit circle, which amplifies the factor $\sqrt{(1+|\alpha_k|)/(1-|\alpha_k|)}$ in \eqref{errbound}, increases the coefficients $\beta_{k,\ell}$ in \eqref{betaklAp} and \eqref{betaklB}, and makes them decay more slowly with $k$ and $\ell$. The resulting error bound is therefore larger and exhibits slower decay, yielding an artificially pessimistic estimate that can substantially underestimate the actual convergence rate associated with the choice of minimal $E$. 
} 

{Meanwhile, this observation also suggests that the sketched algorithm \emph{may} converge more slowly than the non-sketched algorithm if $W(T_m) \not \subset W(A)$. In this situation, the minimal $E$ containing both $W(A)$ and $W(T_m)$ due to sketching is a superset of the minimal $E$ containing $W(A)$ alone in the non-sketched case, since $H_m = V_m^TAV_m$ with orthonormal basis satisfies $W(H_m) \subset W(A)$. A larger minimal $E$ is linked to slower convergence rate. However, as we show in numerical experiments, the actual convergence rate of the non-sketched and the sketched tangential RKSM is largely comparable for a variety of test matrices.} {This is not unexpected for FD-type a priori bounds: the set \(E\) and the pole-separation assumptions are chosen to certify a worst-case convergence rate, and the resulting estimate can be much more conservative than the observed convergence of the projected method.}

\begin{figure}[hbtp!]
\centering
\setlength{\tabcolsep}{2pt} %
\renewcommand{\arraystretch}{0} %

\begin{tabular}{@{}c c@{}}
  \includegraphics[width=.48\textwidth]{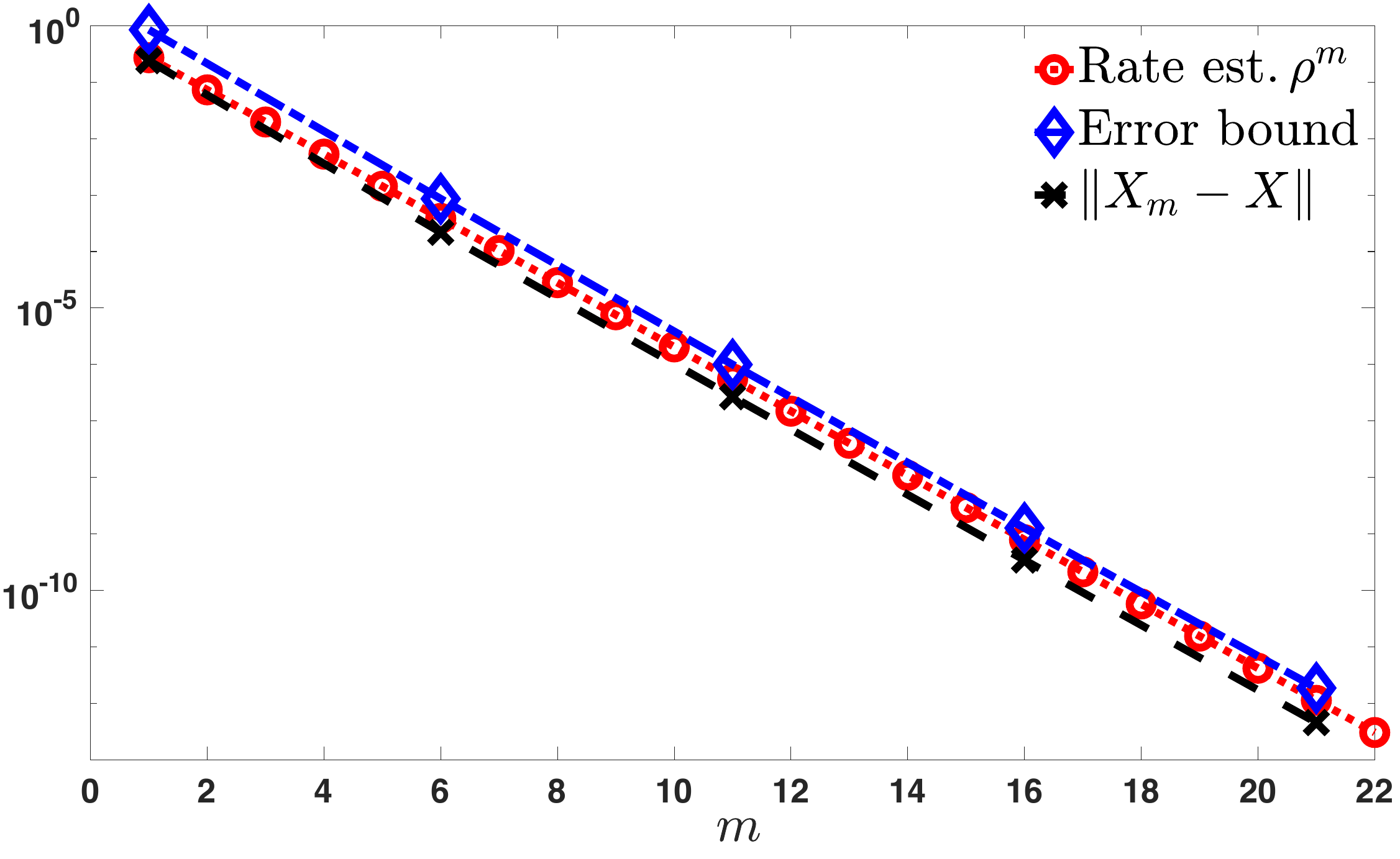} &
  \includegraphics[width=.48\textwidth]{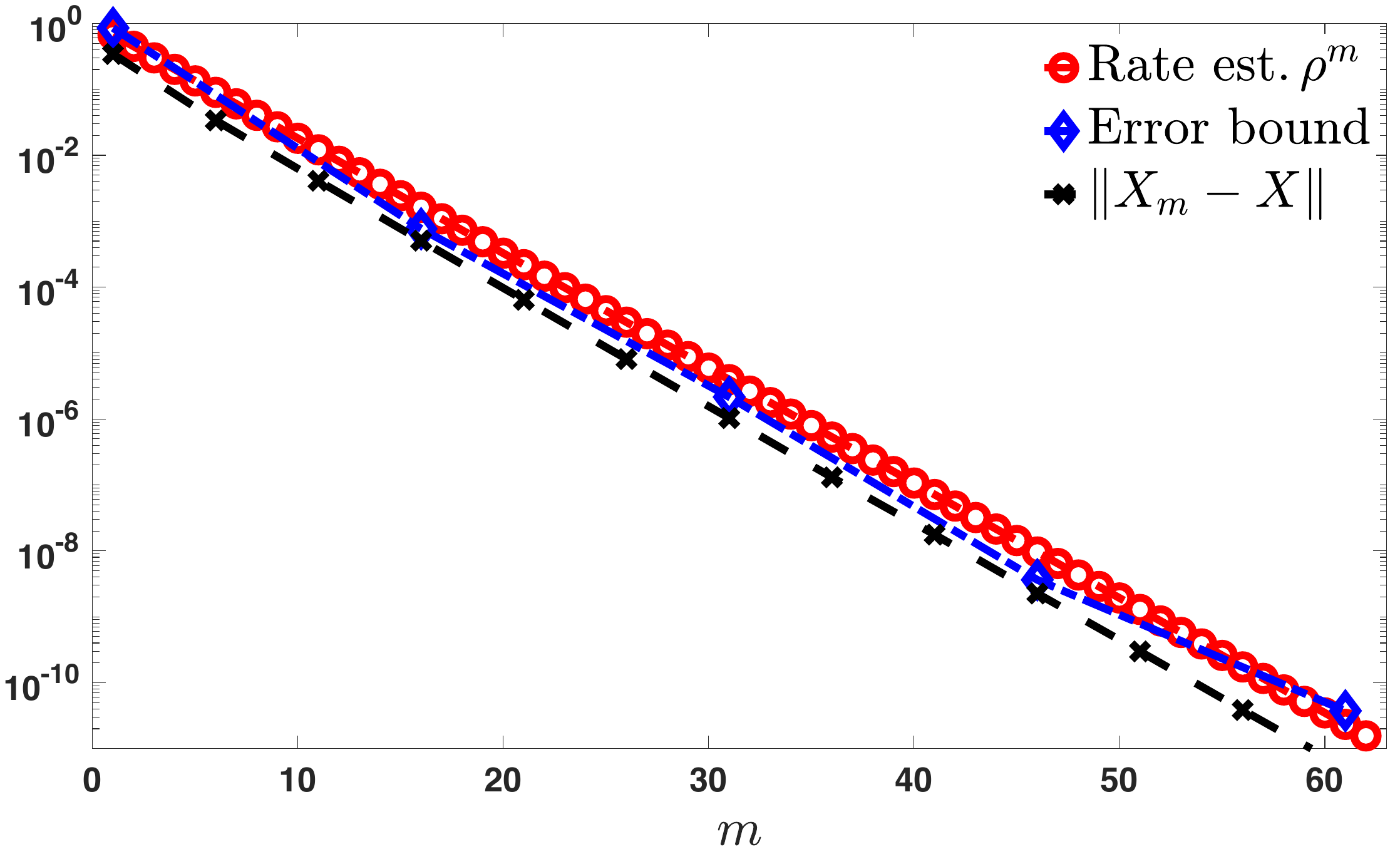} \\
  \includegraphics[width=.48\textwidth]{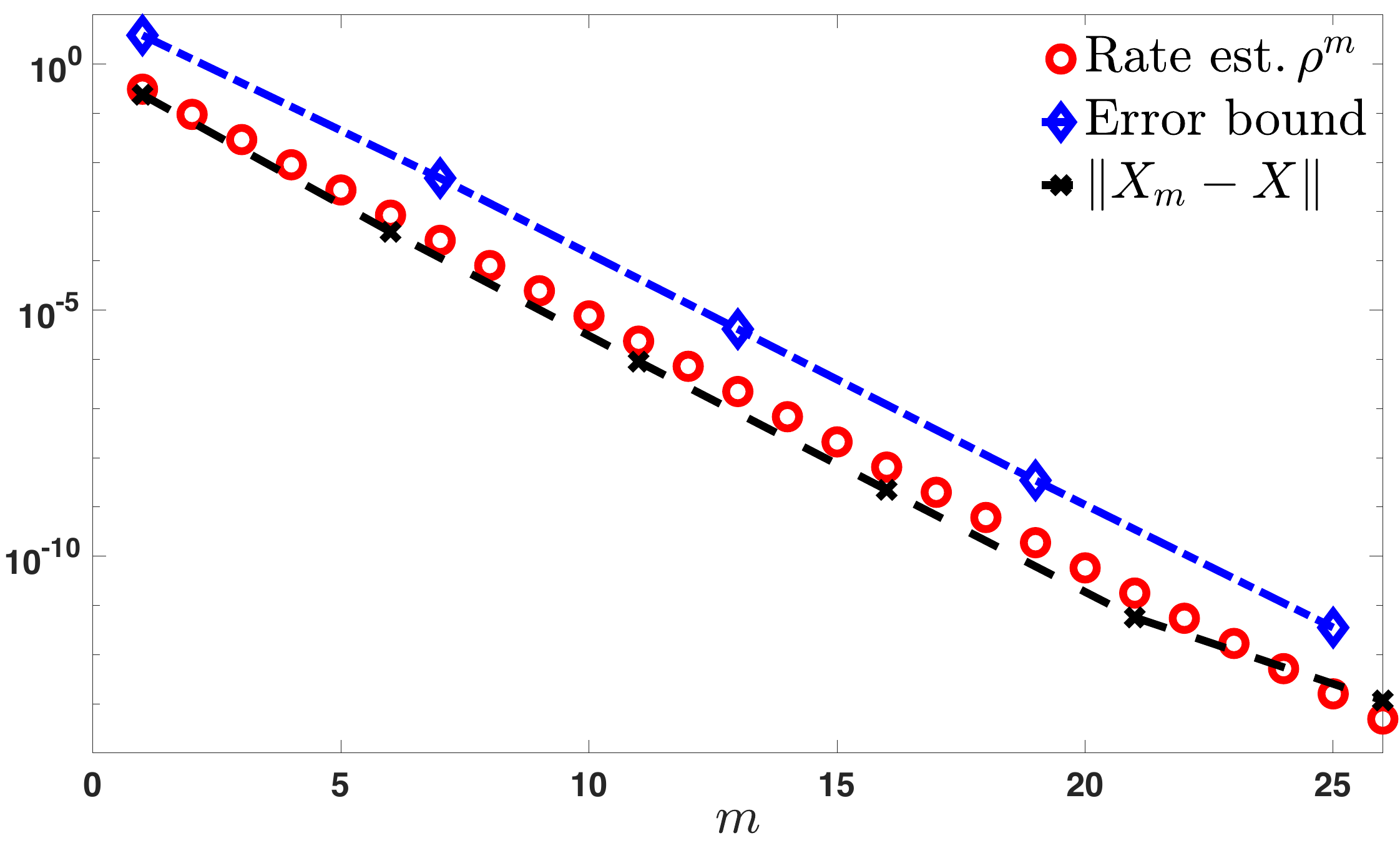} &
  \includegraphics[width=.48\textwidth]{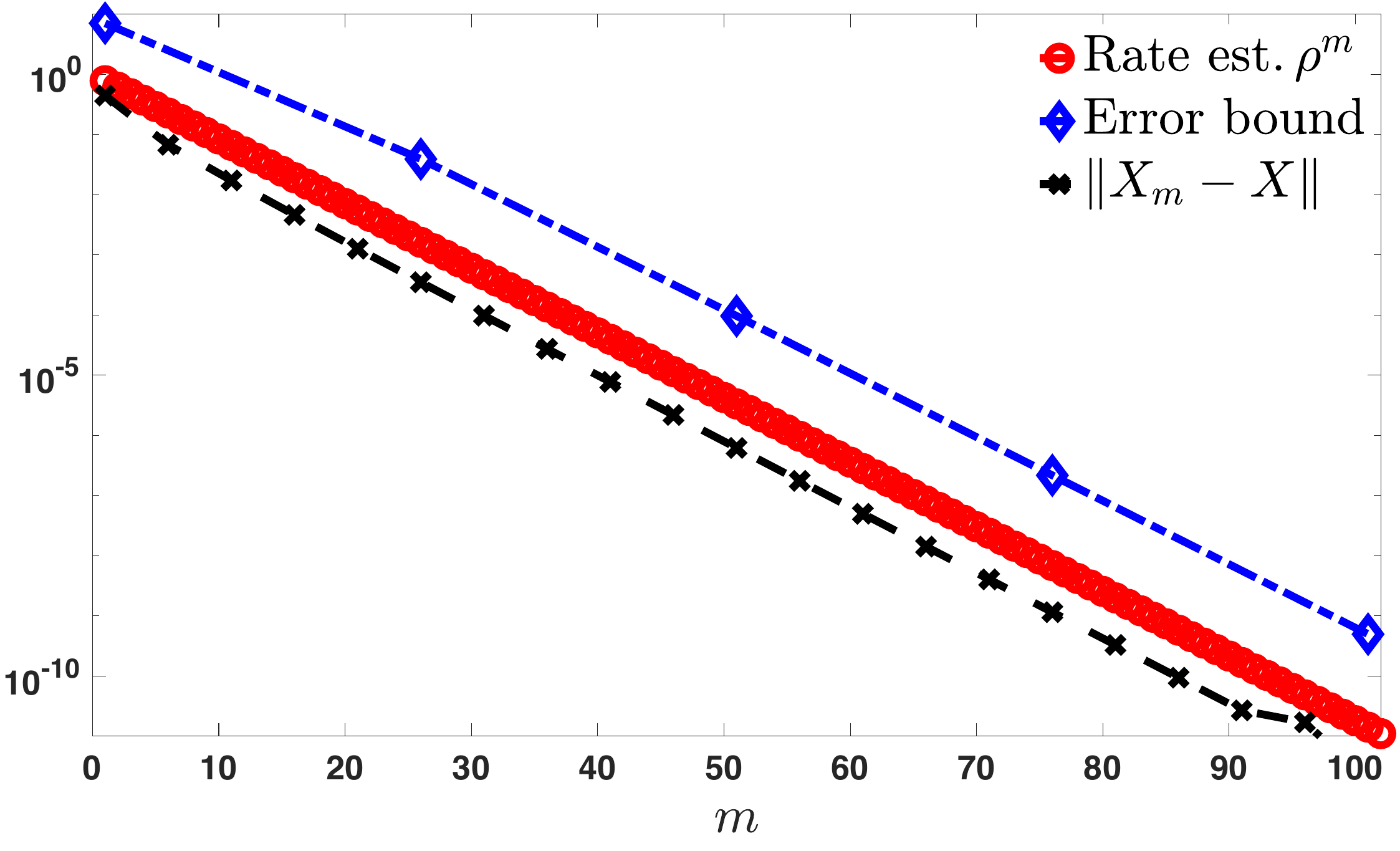}
\end{tabular}

\caption{{Actual errors $\|X_m-X\|$ (black), error bounds \eqref{errbound} on $\|X_m-X\|$ \textup{(}blue\textup{)}, and the error rate estimate $\rho^m$ \textup{(}red\textup{)}. \textup{(}Top left\,:\,$W(A)=[-10^2,-1]$, $\rho \approx 0.26987$; top right\,:\,$W(A)=[-10^4,-1]$, $\rho \approx 0.66942$; bottom left\,:\,$W(A)=\{z:|z+5|\leq 4\}$, $\rho \approx 0.30769$; bottom right\,:\,$W(A) = \{z:|z+30| \leq 29\}$, $\rho \approx 0.78087$\textup{)}.}}
\label{fig:errorbnd}\vspace{-0.2in}
\end{figure}

We give two examples showing that our bound \eqref{errbound} accurately reflects the convergence rate of RKSM. Consider EKSM for solving the Lyapunov equation, where $W(A)$ is (i) a line segment on the negative real axis or (ii) a disk on the left half complex plane. For an $m$-step EKSM that develops a subspace of dimension $2m+1$, the error is bounded by $\mathcal{O}(\rho^m)$, where $\rho$ is the convergence factor. For problem (i), $\rho = \Big(\frac{\kappa^{1/4}-1}{\kappa^{1/4}+1}\Big)^2$, where $\kappa = \frac{\lambda_{max}(-A)}{\lambda_{min}(-A)}$ is the $2$-condition number of $A$; for problem (ii), $\rho = \frac{r^2}{4c^2-3r^2}$, where $c<0$ is the center of disk and $r$ is the radius ($r < |c|$); see \cite{Knizhnerman.Simoncini.2011}. For each problem, we invoke FFT with $2^pN = 2^{14}$ equispaced nodes on $[0,2\pi)$ to evaluate $\beta_{k,\ell}$'s, where $m \leq \ell \leq m+2^{11}$ for all $k$ ($1 \leq k \leq m+2^{11}$) ($\beta_{k,\ell}$ is disregarded for $k,\ell > m+2^{11}$). We plot the bound \eqref{errbound} and the error estimate $\rho^m$ in the semilog scale in Figure~\ref{fig:errorbnd}, showing that the bound correctly reflects the convergence of EKSM, as the two curves are mostly parallel. For more difficult problems with $\rho$ closer to $1$, more equispaced quadrature nodes are needed to evaluate $\beta_{k,\ell}$'s accurately, and the computation of the bound becomes more expensive. {The errors $\|X_m-X\|$ computed by EKSM are also plotted (in black), where matrix $A$ is a $16384\times 16384$ real diagonal or block diagonal (with $2\times 2$ diagonal blocks) sparse matrix with eigenvalues uniformly distributed (in equilibrium measure) in $W(A)$, and $b$ is a normalized vector generated by MATLAB's \texttt{randn} function. The error curves are  consistently parallel to and below the curves of the error bound \eqref{errbound}; they seem to differ by a constant factor that may depend on $W(A)$.}

\section{Numerical Experiments} We carry out an extensive set of numerical experiments comparing the new methods with a few existing methods in the number of iterations, storage cost, and runtime. Note that such comparisons may not be straightforward due to distinct stopping criterion adopted by different methods and the choice of different parameters. Nevertheless, we want to share some experiences about the strength and weakness of different algorithms that are unlikely to be very sensitive to the variation in stopping criterion and parameters. 

\subsection{Test problems} We tested 13 matrices from diverse applications, including model reductions, complex network modeling, and partial differential equations (PDEs). Most matrices can be found by name on SuiteSparse \cite{kolodziej2019suitesparse}, whereas \textit{aerofoil} is a test problem in \cite{Timme.etal.2012},  provided by Sebastian Timme at the University of Liverpool.  \textit{Laplace2D} is a discrete Laplacian on a unit square based on the standard 5-point stencil finite difference discretization, with $512$ uniform subintervals in each direction and with mesh size $h = \frac{1}{512}$. \textit{Laplace3D} is a discrete Laplacian on a unit cube based on the standard 7-point stencil finite difference discretization, with $65$ uniform subintervals in each direction and with mesh size $h = \frac{1}{65}$. %
For the two problems of undirected graphs, we let $A$ be the graph Laplacian of their adjacency matrices, without scaling. Since the graph Laplacians have 1-dimensional null space spanned by ${\bf 1} = [1,1,\ldots,1]^T$, the original right-hand side vectors $B$ in the Lyapunov equation have been orthogonalized against this null space basis vector to ensure a unique solution. The symmetry, size, number of nonzero elements (nnz), and the estimate of the intersection of the numerical range $W(-A)$ and the real axis of these matrices are given in Table~\ref{problemdscp}. For the sketched method, we used the sparse dimension reduction map with $\epsilon = 10^{-0.5} \approx 0.3162$ in the sketching assumption \eqref{sketchassumption}. For a subspace $\mathrm{span}(U_m) \subset \mathbb{R}^n$ of dimension $m \ll n$, the mapped dimension is $s = \frac{m}{\epsilon^2} = 10m$. 

All the experiments were performed on an iMac Desktop with 3.2 GHz Quad-Core Intel Core i5 processor, and 32 GB 1867 MHz DDR3 memory, running MacOS Monterey (12.7.6) and MATLAB R2023b Update 8. 

\begin{table}[hbtp!]
{\footnotesize
   \begin{center}
  \begin{tabular}{|r|r|c|r|r|r|r|}
    \hline
    problem & application & symmetric & size & nnz & spectral interval est.\\ \hline
    \textit{gyro\_m} & model reduction & Yes & $17361$ & $340431$ & $[7\times 10^{-14}, 2\times 10^{-7}]$ \\ 
    \textit{chipcool1} & model reduction & No & $20082$ & $281150$ & $[2\times 10^{-7},0.32]$ \\ 
    \textit{rail\_20209} & model reduction & Yes & $20209$ & $139233$ & $[2\times10^{-10},6\times 10^{-5}]$ \\ 
    \textit{rail\_79841} & model reduction & Yes & $79841$ & $553921$ & $[5\times10^{-3},7\times10^{3}]$ \\
    \textit{gas\_sensor} & model reduction & Yes & $66917$ & $1703365$ & $[9\times10^{-8},21]$ \\
    \textit{filter3D} & model reduction & Yes & $106437$ & $2707179$ & $[1.5\times 10^{-9},1.6\times 10^{-2}]$ \\ 
    \textit{boneS01} & model reduction & Yes & $127224$ & $5516602$ & $[2\times10{-3},5\times10^{4}]$ \\ 
    \textit{wing} & undirected graph & Yes & $62032$ & $243088$ & $[10^{-3},8]$ \\ 
    \textit{delaunay\_n16} & undirected graph & Yes & $65536$ & $393150$ & $[4\times10^{-4},19]$ \\
    \textit{aerofoil} & aerodynamics & No & $16388$ & $627979$ & $[3\times 10^{-3},2\times 10^{3}]$ \\
    \textit{af23560} & fluid dynamics & No & $23560$ & $460598$ & $[0.2,2\times 10^{2}]$ \\
    \textit{Laplace2D} & PDEs & Yes & $261121$ & $1303561$ & $[2\pi^2,2.1 \times 10^{6}]$ \\
    \textit{Laplace3D} & PDEs & Yes & $262144$ & $1810432$ & $[3\pi^2,5.1\times10^{4}]$ \\ \hline
    \end{tabular}
    \end{center}
  }
  \caption{Test Matrices for Numerical Experiments}\label{problemdscp}\vspace{-0.25in}
\end{table}

\subsection{Solving Lyapunov equations with a low-rank \({B}\)} Let us first consider solving the Lyapunov equations $AX+XA^T+BB^T=0$, where $B$ has a low column rank that is typically no greater than $5$. For such an equation, tangential RKSM usually does not need many iteration steps or a large subspace to converge. As a result, the orthogonalization of basis vectors is usually not the dominant cost of the algorithm. This is particularly the case if the linear solves performed at each RKSM step are expensive. In this circumstance, if the linear solves can still be performed efficiently by a sparse direct solver, our sketched tangential RKSM with a few repeated poles tends to have an advantage in runtime over the standard tangential RKSM. 

We chose five problems where the linear solves involving the shifted matrices for RKSM can be performed by MATLAB's direct solvers at a relatively high cost. We computed the sparse LU or Cholesky factorization of the matrix $sI-A$, where $s$ is the lower end of the estimated real interval covered by $W(-A)$, and we evaluated the density of the sparse factorizations, defined as $\nu = \frac{nnz(L)+nnz(U)}{n}$ (by MATLAB's \texttt{lu} function with five output variables, if $A$ is non-symmetirc) or $\nu = \frac{2nnz(L)}{n}$ (by the \texttt{chol} function with symmetric approximate minimum degree permutation, if $A$ is symmetric), shown in Table~\ref{test1}. This quantity is the average number of nonzeros per row in the application of the triangular linear solves. A higher density means that performing the sparse direct solver is more expensive relative to the problem size. 

For the right-hand side block $B$, we used the built-in data for problems \textit{chipcool1} and \textit{gas\_sensor}. These blocks are referred to as matrix $C$, which contain 5 columns for \textit{chipcool1} and 28 columns for \textit{gas\_sensor}, respectively. {For \textit{boneS01}, the matrix $B = \widehat{E}_p$ is generated by MATLAB's random sampling of $p$ columns of the identity matrix $I_n$ by the \texttt{randsample} function. For \textit{aerofoil}, $B = \widehat{I}_p = [I_p; \, 0_{(n-p)\times p}]$ is the first $p$ columns of the identity matrix $I_n$.} For \textit{Laplace3D}, $B$ contains random elements satisfying the standard normal distribution, generated by MATLAB's \texttt{randn} function. Randomized functions are initiated by the default setting (\texttt{rng default}). 

\begin{table}[hbtp!]
{\footnotesize
   \begin{center}
  \begin{tabular}{|r|r|r|r|r|r|r|r|r|r|}
    \hline
    problem & $B$ & $m$ & \multicolumn{3}{c|}{Standard TRKSM} & \multicolumn{4}{c|}{sketched 6-pole TRKSM} \\ \cline{4-10}
    \& density &  & & step & dim & runtime & step & dim & \multicolumn{2}{c|}{runtime \& ratio}\\ \hline
    
    \textit{chipcool1} & $C(:,1)$ & $30$ & $19$ & $20$ & $9.69$ & $23$ & $24$ & $5.40$ & $55.7\%$ \\ 
    $\nu \approx608$ & $C(:,1:3)$ & $65$ & $20$ & $53$ & $9.97$ & $32$ & $56$ & $5.22$ & $52.4\%$ \\ 
    $n=20082$ & $C(:,1:5)$ & $100$ & $20$ & $83$ & $10.44$ & $38$ & $93$ & $5.70$ & $54.6\%$  \\ \hline
    
    \textit{gas\_sensor} & $C(:,1)$ & $45$ & $29$ & $30$ & $48.90$ & $44$ & $45$ & $19.12$ & $39.1\%$ \\
    $\nu \approx1674$ & $C(:,1:3)$ & $130$ & $29$ & $83$ & $51.12$ & $56$ & $110$ & $21.73$ & $42.5\%$ \\
    $n=66917$ & $C(:,1:5)$ & $200$ & $27$ & $129$ & $51.64$ & $61$ & $178$ & $27.13$ & $52.5\%$  \\ \hline
   
    \textit{boneS01} & $\widehat{E}_2$ & $45$ & $21$ & $43$ & $43.02$ & $22$ & $44$ & $17.06$ & $39.7\%$ \\ 
    $\nu \approx587$ & $\widehat{E}_4$ & $90$ & $20$ & $82$ & $43.84$ & $24$ & $92$ & $18.72$ & $42.7\%$ \\ 
    $n=127224$ & $\widehat{E}_8$ & $190$ & $18$ & $150$ & $50.86$ & $25$ & $188$ & $26.75$ & $52.6\%$ \\ \hline
    
    \textit{aerofoil} & $\widehat{I}_1$ & $50$ & $48$ & $49$ & $11.03$ & $47$ & $48$ & $2.66$&  $24.1\%$ \\
    $\nu \approx402$ & $\widehat{I}_2$ & $90$ & $48$ & $88$ & $11.83$ & $58$ & $91$ & $3.36$ & $28.4\%$ \\
    $n=16388$ & $\widehat{I}_3$ & $150$ & $50$ & $132$ & $13.50$ & $63$ & $142$ & $4.01$ & $29.7\%$ \\ \hline
    
    \textit{Laplace3D} & $\mathrm{randn}(n,2)$ & $30$ & $12$ & $26$ & $95.35$ & $13$ & $28$ & $31.35$ & $32.9\%$ \\
    $\nu \approx1405$ & $\mathrm{randn}(n,4)$ & $60$ & $13$ & $56$ & $105.68$ & $12$ & $52$ & $33.01$ & $31.2\%$ \\
    $n=262144$ & $\mathrm{randn}(n,8)$ & $120$ & $12$ & $104$ & $106.63$ & $12$ & $104$ & $40.22$ & $37.7\%$  \\ \hline
    \end{tabular}
    \end{center}
  }
  \caption{Comparison of the standard and the sketched multi-directional TRKSM with six repeated poles for solving $AX+XA^T+BB^T=0$ with a low-rank $B$. {Matrix $C$ is a built-in matrix of the problem data; $\widehat{E}_p$ contains random sampling of $p$ columns of the identity matrix $I_n$ by MATLAB's \texttt{randsample} function; $\widehat{I}_p = [I_p; \, 0_{(n-p)\times p}]$ is the first $p$ columns of the identity matrix $I_n$.}}\label{test1} \vspace{-0.25in}
\end{table}

We let the tolerance of the relative residual norm be $10^{-6}$. The largest rational Krylov subspace dimension $m$ chosen for each problem is shown in Table~\ref{test1}. The depth of selective orthogonalization for the sketched method is $\frac{m}{4}=0.25m$. One can see from the table that the sketched tangential RKSM using six repeated poles outperforms the standard algorithm in runtime, though the former takes more steps and therefore needs a larger subspace to converge. The six poles were determined in the first six RKSM steps adaptively, described in line 9 of Algorithm~\ref{alg:rstgrksm}. The last column in Table~\ref{test1}, marked as runtime \& ratio, is the runtime of the sketched method as a percentage of the runtime of the standard method. The advantage in runtime is achieved by reusing the existing sparse LU or Cholesky factorizations obtained in the first six RKSM steps. The improvement is also attributed to the small subspaces, which do not incur a high cost of orthogonalization. 

Meanwhile, it should be noted that the advantage in runtime of the 6-pole algorithm becomes less pronounced for $B$ of higher rank. As the rank of $B$ increases, the standard adaptive tangential RKSM takes roughly the same number of iteration steps (i.e., number of direct linear solves with multiple right-hand sides) to converge, whereas the 6-pole sketched variant needs more steps and a larger subspace. This entails higher cost of linear solves, orthogonalization, and storage for the latter. To handle $B$ of higher rank (e.g., 10--15 and above), we would use the sketched algorithm with fully adaptive poles, without saving sparse matrix factorizations. 

\begin{figure}[htbp!]
\centering
\begin{tikzpicture}

\def\cOrth{blue!60}
\def\cSolve{orange!85}
\def\cRes{green!60!black!30}
\def\cPole{violet!60}
\def\cSketch{red!70}
\def\cOther{gray!35}

\pgfplotsset{
  runtimeplotfive/.style={
    ybar stacked,
    bar width=7pt,
    enlarge x limits=0,
    width=0.55\linewidth,
    height=8.5cm,
    ymin=0,
    ymax=240,
    enlarge y limits=false,
    ylabel={Runtime (sec)},
    ymajorgrids=true,
    grid style={dashed},
    legend style={
      font=\scriptsize,
      at={(0.5,-0.20)},
      anchor=north,
      legend columns=3,
      draw=none,
      legend cell align={left}
    },
    symbolic x coords={
      LeftGap,
      chipStd,chipSk,chipGap,
      gasStd,gasSk,gasGap,
      boneStd,boneSk,boneGap,
      aeroStd,aeroSk,aeroGap,
      L3Std,L3Sk,L3Gap
    },
    xtick={chipStd,gasStd,boneStd,aeroStd,L3Std},
    xticklabels={
      \texttt{chipcool1},
      \texttt{gas\_sensor},
      \texttt{boneS01},
      \texttt{aerofoil},
      \texttt{Lpalace3D}
    },
    xticklabel style={font=\scriptsize, rotate=40, anchor=east},
  }
}

\begin{axis}[runtimeplotfive]

\addplot[fill=\cOrth, draw=none] coordinates {
(LeftGap,0)
(chipStd,0.512) (chipSk,0.109) (chipGap,0)
(gasStd,5.314) (gasSk,0.860) (gasGap,0)
(boneStd,12.844) (boneSk,1.332) (boneGap,0)
(aeroStd,1.367) (aeroSk,0.218) (aeroGap,0)
(L3Std,15.195) (L3Sk,1.519) (L3Gap,0)
};
\addlegendentry{Basis Orth}

\addplot[fill=\cSolve, draw=none] coordinates {
(LeftGap,0)
(chipStd,20.252) (chipSk,10.097) (chipGap,0)
(gasStd,111.364) (gasSk,58.611) (gasGap,0)
(boneStd,67.365) (boneSk,53.364) (boneGap,0)
(aeroStd,23.857) (aeroSk,5.401) (aeroGap,0)
(L3Std,207.564) (L3Sk,104.169) (L3Gap,0)
};
\addlegendentry{Lin Solve}

\addplot[fill=\cRes, draw=none] coordinates {
(LeftGap,0)
(chipStd,0.527) (chipSk,0.077) (chipGap,0)
(gasStd,3.005) (gasSk,0.515) (gasGap,0)
(boneStd,4.956) (boneSk,0.224) (boneGap,0)
(aeroStd,2.494) (aeroSk,0.262) (aeroGap,0)
(L3Std,5.108) (L3Sk,0.109) (L3Gap,0)
};
\addlegendentry{Res Estimate}

\addplot[fill=\cPole, draw=none] coordinates {
(LeftGap,0)
(chipStd,0.338) (chipSk,1.041) (chipGap,0)
(gasStd,0.714) (gasSk,2.813) (gasGap,0)
(boneStd,0.521) (boneSk,1.037) (boneGap,0)
(aeroStd,1.927) (aeroSk,1.769) (aeroGap,0)
(L3Std,0.285) (L3Sk,0.406) (L3Gap,0)
};
\addlegendentry{New Pole}

\addplot[fill=\cSketch, draw=none] coordinates {
(LeftGap,0)
(chipStd,0) (chipSk,0.213) (chipGap,0)
(gasStd,0) (gasSk,1.166) (gasGap,0)
(boneStd,0) (boneSk,1.751) (boneGap,0)
(aeroStd,0) (aeroSk,0.414) (aeroGap,0)
(L3Std,0) (L3Sk,1.806) (L3Gap,0)
};
\addlegendentry{Sketching}

\addplot[fill=\cOther, draw=none] coordinates {
(LeftGap,0)
(chipStd,0.014) (chipSk,0.104) (chipGap,0)
(gasStd,0.050) (gasSk,0.611) (gasGap,0)
(boneStd,0.086) (boneSk,0.245) (boneGap,0)
(aeroStd,0.171) (aeroSk,0.323) (aeroGap,0)
(L3Std,0.083) (L3Sk,0.248) (L3Gap,0)
};
\addlegendentry{Other Costs}

\node[font=\scriptsize, anchor=west] at (rel axis cs:0.01,0.98)
{left bar = TRKSM, right bar = S-TRKSM};

\end{axis}
\end{tikzpicture}
\caption{Runtime breakdown (stacked) for TRKSM (left bar) and S-TRKSM (right bar) on 5 test problems, with matrix $B$ as the last choice of the three in Table~\ref{test1}. Other cost accounts for unreported overhead so each stack matches the total runtime.}
\label{fig:runtime_stacked_5problems}\vspace{-0.2in}
\end{figure}
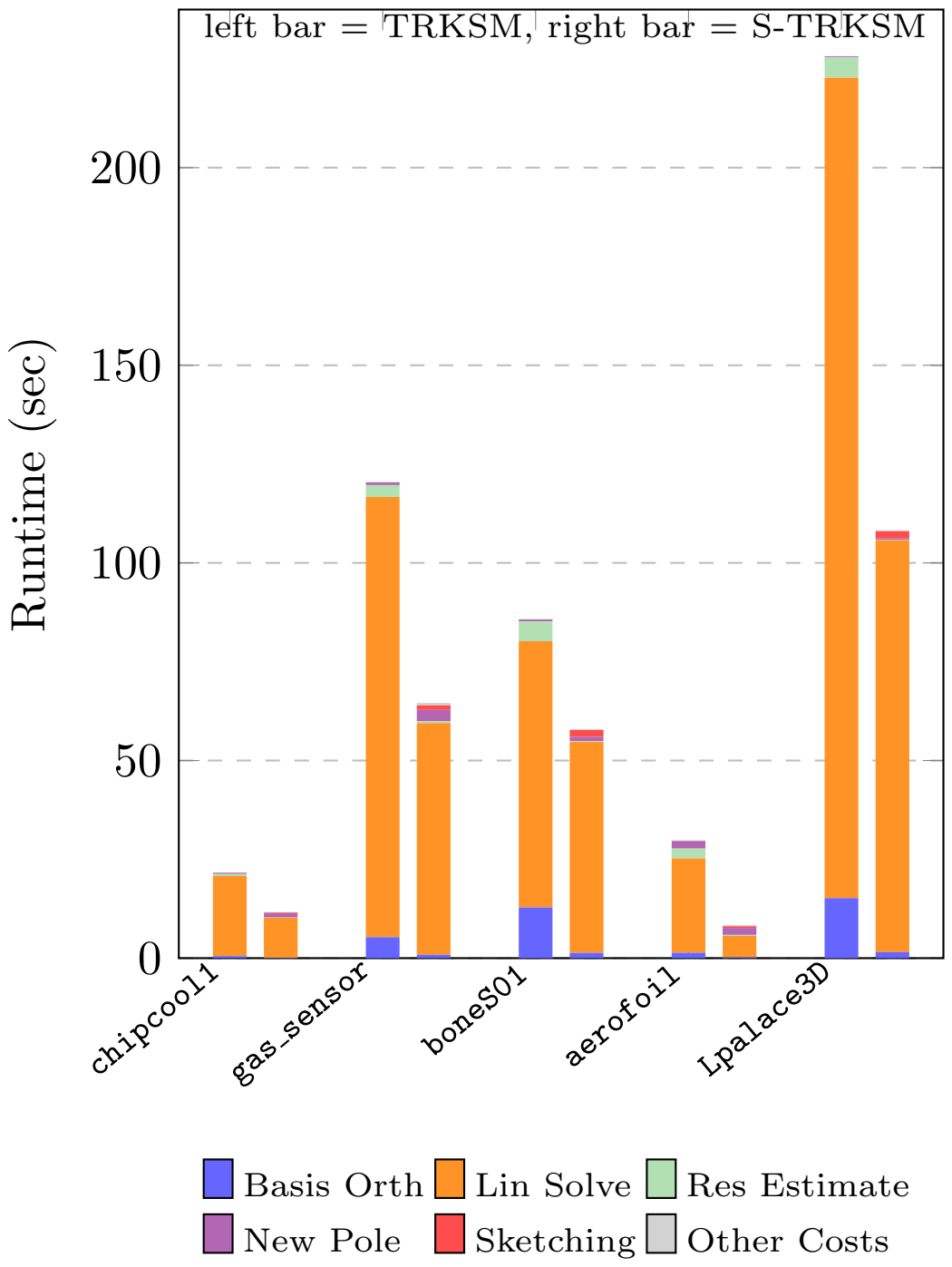

{In Figure~\ref{fig:runtime_stacked_5problems}, a stacked bar plot is drawn to show the details of the runtime of different computational tasks. For each problem, we let the matrix $B$ be the last of the three choices in Table~\ref{test1}; the left and the right stacked bars represent the runtime of the standard and the sketched algorithms, respectively. The runtime recorded include that for basis orthogonalization (blue), linear solves (orange), equation residual norm estimate (light green), search for new poles (violet), sketching-related overheads (red), and other costs (gray). Since the maximum subspace dimension used is not very high, the total runtime is dominated by the linear solves. By reusing the LU or Cholesky factorizations obtained in the first few iterations, the runtime is reduced effectively. The change in other colored bars, and basis orthogonalization in particular, will be seen more clearly with a higher-rank $B$ shown in the next subsection.}

Another interesting observation here is that the use of a small number of repeated poles may generate approximate solutions of higher accuracy. In Table~\ref{test2}, we show a few examples where the 3-pole or 4-pole sketched algorithm is able to achieve a considerably smaller relative residual norm than the standard variant based on adaptive poles. The residual norms shown are the lowest achievable level for each algorithm. We also ran the sketched variant with adaptive poles and saw stagnation in the residual norm similar to that of the standard algorithm with adaptive poles. We found that in the sequence of adaptive poles, many are small ones (close to the origin), whereas in the sequence of repeated poles, small poles are used only two or three times in the early stage of RKSM steps, followed by poles that are considerably larger. The shifted matrices $sI-A$ with smaller poles $s$ have a higher condition number than those with larger poles, especially if $A$ is ill-conditioned. This is probably a contributing factor of the difference in the achievable residual norm. 

\begin{table}[hbtp!]
{\footnotesize
   \begin{center}
  \begin{tabular}{|r|r|r|r|r|r|r|r|r|r|}
    \hline
    problem  & $B$ & $m$ & \multicolumn{3}{c|}{MultiDir TRKSM} & \multicolumn{4}{c|}{sketched MultiDir TRKSM} \\ %
    &  & & \multicolumn{3}{c|}{(with adaptive poles)} & \multicolumn{4}{c|}{(with a few repeated poles)} \\ \cline{4-10}
    & & & step & dim & residual & poles & step & dim & residual \\ \hline
    
    \textit{aerofoil} & $\widehat{I}_2$ & $250$ & $51$ & $94$ & $4.8\times 10^{-7}$ & $4$ & $158$ & $236$ & $4.5\times 10^{-9}$ \\ \hline
    
    \textit{af23560} & $\widehat{I}_2$ & $125$ & $40$ & $82$ & $5.2\times 10^{-7}$ & $3$ & $53$ & $108$ & $8.5 \times 10^{-11}$ \\
     & $\widehat{E}_2$ & $160$ & $49$ & $100$ & $4.1\times 10^{-5}$ & $3$ & $76$ & $150$ & $8.9 \times 10^{-11}$ \\ \hline
    \end{tabular}
    \end{center}
  }
  \caption{Examples of improved accuracy of approximate solutions computed by the sketched TRKSM with a few repeated poles}\label{test2}\vspace{-0.25in}
\end{table}

\subsection{Solving Lyapunov equations with $B$ of higher rank} For the right-hand side block $B$ of higher rank (typically $15$ or higher), we compare the standard and the sketched tangential RKSM with adaptive poles. The tolerance of the relative residual norm is $10^{-7}$ for all problems from model reduction and undirected graphs, and it is $5\times 10^{-6}$ and $2 \times 10^{-6}$ for \textit{aerofoil} and \textit{af23560}, respectively. All linear solves are performed by MATLAB's backslash operator, even though this is not necessarily the optimal solver for each problem. We keep this choice to illustrate the impact of the cost of linear solvers on the runtime improvement of the sketched algorithm. 

\begin{table}
{\footnotesize
   \begin{center}
  \begin{tabular}{|r|r|r|r|r|r|r|r|r|}
    \hline
    problem size & $B$ & \multicolumn{3}{c|}{MultiDir TRKSM} & \multicolumn{4}{c|}{sketched MultiDir TRKSM} \\ \cline{3-9}
    \& density &  & step & dim & runtime & step & dim & \multicolumn{2}{c|}{runtime \& ratio}\\ \hline
    \textit{gyro\_m} & $\mathrm{randn}(n,15)$ & $30$ & $465$ & $21.18$ & $35$ & $468$ & $14.20$ & $67.0\%$ \\ 
    $n=17361$ & $\mathrm{randn}(n,30)$ & $28$ & $869$ & $91.32$ & $29$ & $900$ & $55.66$ & $61.0\%$ \\ 
    $\nu \approx48$ & $A^{-1}\widehat{I}_{90}$ & $10$ & $422$ & $16.99$ & $12$ & $479$ & $11.32$ & $66.6\%$ \\ 
    & $A^{-1}\widehat{I}_{300}$ & FAIL & $-$ & $-$ & $10$ & $867$ & $42.55$ & $-$ \\ \hline
    \textit{chipcool1} & $[C\:\widehat{E}_{20}]$ & $26$ & $450$ & $35.33$ & $38$ & $458$ & $34.63$ & $98.0\%$ \\ 
    $n=20082$ & $[C\:\widehat{E}_{50}]$ & $26$ & $850$ & $110.72$ & $28$ & $824$ & $70.74$ & $63.9\%$ \\ 
    $\nu \approx608$ & $A^{-1}\widehat{I}_{40}$ & $23$ & $431$ & $36.10$ & $26$ & $488$ & $31.46$ & $87.1\%$ \\ 
    & $A^{-1}\widehat{I}_{80}$ & $24$ & $859$ & $110.44$ & $24$ & $838$ & $71.44$ & $64.7\%$ \\ \hline
    \textit{rail\_20209} & $\mathrm{randn}(n,30)$ & $15$ & $480$ & $20.99$ & $15$ & $476$ & $9.32$ & $44.4\%$ \\ 
    $n=20209$ & $\mathrm{randn}(n,60)$ & $14$ & $884$ & $68.67$ & $14$ & $881$ & $34.00$ & $49.5\%$ \\
    $\nu \approx33$ & $A^{-1}\widehat{I}_{65}$ & $8$ & $449$ & $15.41$ & $9$ & $507$ & $9.45$ & $61.3\%$ \\ 
    & $A^{-1}\widehat{I}_{140}$ & $8$ & $917$ & $64.33$ & $7$ & $916$ & $35.53$ & $55.2\%$ \\ \hline
    \textit{rail\_79841} & $\mathrm{randn}(n,25)$ & $17$ & $442$ & $55.19$ & $17$ & $442$ & $13.88$ & $25.1\%$ \\
    $n=79841$ & $\mathrm{randn}(n,55)$ & $16$ & $909$ & $228.18$ & $17$ & $960$ & $62.27$ & $27.3\%$ \\
    $\nu \approx43$ & $A^{-1}\widehat{I}_{60}$ & $9$ & $485$ & $57.87$ & $10$ & $532$ & $14.39$ & $24.9\%$ \\
    & $A^{-1}\widehat{I}_{130}$ & $9$ & $1013$ & $232.94$ & $8$ & $1005$ & $52.79$ & $22.7\%$ \\ \hline
    \textit{gas\_sensor} & $C(:,1:14)$ & $31$ & $398$ & $102.53$ & $49$ & $432$ & $100.81$ & $98.3\%$ \\
    $n=66917$ & $C$ & $32$ & $790$ & $215.33$ & $51$ & $795$ & $161.88$ & $75.2\%$ \\
    $\nu \approx1674$ & $A^{-1}\widehat{I}_{20}$ & $28$ & $514$ & $116.18$ & $29$ & $508$ & $70.92$ & $61.0\%$ \\ 
    & $A^{-1}\widehat{I}_{40}$ & $28$ & $969$ & $290.59$ & $29$ & $982$ & $127.77$ & $44.0\%$ \\ \hline
    \textit{filter3D} & $\mathrm{randn}(n,20)$ & $22$ & $458$ & $95.68$ & $22$ & $456$ & $34.36$ & $35.9\%$ \\ 
    $n=106437$ & $\mathrm{randn}(n,45)$ & $20$ & $941$ & $336.50$ & $20$ & $940$ & $85.73$ & $25.5\%$ \\ 
    $\nu \approx342$ & $A^{-1}\widehat{I}_{45}$ & $22$ & $454$ & $99.72$ & $23$ & $453$ & $37.32$ & $37.4\%$ \\ 
    & $A^{-1}\widehat{I}_{90}$ & $22$ & $897$ & $323.49$ & $22$ & $893$ & $82.67$ & $25.6\%$ \\ \hline
    \textit{boneS01} & $\widehat{E}_{25}$ & $19$ & $485$ & $134.52$ & $19$ & $483$ & $55.57$ & $41.3\%$ \\ 
    $n=127224$ & $\widehat{E}_{50}$ & $17$ & $863$ & $359.74$ & $18$ & $883$ & $100.35$ & $27.9\%$ \\ 
    $\nu \approx587$ & $A^{-1}\widehat{I}_{100}$ & $10$ & $579$ & $145.52$ & $8$ & $462$ & $30.69$ & $21.1\%$ \\
    & $A^{-1}\widehat{I}_{400}$ & $7$ & $985$ & $346.03$ & $6$ & $960$ & $66.98$ & $19.4\%$ \\ \hline

    \textit{wing} & $\mathrm{randn}(n,30)$ & $13$ & $420$ & $38.94$ & $13$ & $420$ & $12.99$ & $33.4\%$ \\ 
    $n=62032$ & $\mathrm{randn}(n,65)$ & $13$ & $910$ & $167.63$ & $13$ & $910$ & $49.65$ & $29.6\%$ \\
    $\nu  \approx 70$ & $A^{-1}\widehat{I}_{50}$ & $12$ & $498$ & $51.51$ & $11$ & $450$ & $12.73$ & $24.7\%$ \\
    & $A^{-1}\widehat{I}_{100}$ & $11$ & $894$ & $156.87$ & $12$ & $993$ & $53.98$ & $34.4\%$ \\ \hline
    \textit{delaunay\_n16} & $\mathrm{randn}(n,30)$ & $15$ & $478$ & $50.16$ & $16$ & $508$ & $15.47$ & $30.8\%$ \\ 
    $n=65536$ & $\mathrm{randn}(n,60)$ & $14$ & $899$ & $179.16$ & $15$ & $960$ & $55.83$ & $31.2\%$ \\ 
    $\nu  \approx 50$ & $A^{-1}\widehat{I}_{35}$ & $13$ & $438$ & $43.49$ & $14$ & $471$ & $13.15$ & $30.2\%$ \\
    & $A^{-1}\widehat{I}_{75}$ & $13$ & $927$ & $185.03$ & $13$ & $917$ & $48.45$ & $26.2\%$ \\ \hline 
    
    \textit{aerofoil} & $\widehat{I}_{20}$ & $55$ & $484$ & $36.04$ & $55$ & $475$ & $27.16$&  $75.4\%$ \\
    $n=16388$ & $\widehat{I}_{40}$ & FAIL & $-$ & $-$ & $58$ & $746$ & $57.94$ & $-$ \\
    $\nu \approx402$ & $A^{-1}\widehat{I}_{25}$ & $53$ & $442$ & $41.24$ & $56$ & $422$ & $26.53$ & $64.3\%$ \\
    & $A^{-1}\widehat{I}_{60}$ & $67$ & $825$ & $113.80$ & $77$ & $842$ & $86.33$ & $75.9\%$ \\ \hline
    \textit{af23560} & $\widehat{I}_{15}$ & FAIL & $-$ & $-$ & $40$ & $448$ & $15.95$ & $-$ \\
    $n=23560$ & $\widehat{I}_{35}$ & FAIL & $-$ & $-$ & $43$ & $903$ & $72.06$ & $-$ \\
    $\nu \approx353$ & $A^{-1}\widehat{I}_{25}$ & $28$ & $512$ & $46.65$ & $28$ & $515$ & $25.64$ & $54.9\%$ \\
    & $A^{-1}\widehat{I}_{50}$ & $28$ & $860$ & $124.66$ & $30$ & $844$ & $52.98$ & $42.5\%$ \\ \hline
    \textit{Laplace2D} & $\mathrm{randn}(n,30)$ & $16$ & $510$ & $216.97$ & $17$ & $511$ & $35.45$ & $16.3\%$ \\
    $n=261121$ & $\mathrm{randn}(n,60)$ & $15$ & $960$ & $706.91$ & $15$ & $960$ & $90.47$ & $12.8\%$ \\
    $\nu \approx72$ & $A^{-1}\widehat{I}_{60}$ & $19$ & $458$ & $196.44$ & $19$ & $458$ & $33.25$ & $16.9\%$ \\
    & $A^{-1}\widehat{I}_{130}$ & $19$ & $958$ & $708.25$ & $19$ & $935$ & $102.51$ & $14.5\%$ \\ \hline
    \textit{Laplace3D} & $\mathrm{randn}(n,35)$ & $11$ & $420$ & $228.55$ & $11$ & $420$ & $105.65$ & $46.2\%$ \\
    $n=262144$ & $\mathrm{randn}(n,70)$ & $11$ & $840$ & $634.88$ & $11$ & $840$ & $162.01$ & $25.5\%$ \\
    $\nu \approx1405$ & $A^{-1}\widehat{I}_{60}$ & $12$ & $494$ & $299.71$ & $12$ & $495$ & $124.04$ & $41.4\%$ \\
    & $A^{-1}\widehat{I}_{120}$ & $12$ & $909$ & $725.28$ & $12$ & $909$ & $180.96$ & $25.0\%$ \\ \hline
    \end{tabular}
    \end{center}
  }
  \caption{Comparison of the standard and the sketched multi-directional TRKSM for solving $AX+XA^T+BB^T=0$ with $B$ of higher rank}\label{test3}
\end{table}

First of all, we underscore the circumstances under which our sketched tangential RKSM is expected to have a decisive advantage in runtime over the standard counterpart. These contributing factors include a larger problem size, a lower cost of linear solves (compared to orthogonalization), and the need for a subspace of higher dimension. For a problem of size $n$, we expect that the sketched vectors belong to $\mathbb{R}^s$ with $s \ll n$, such that the orthogonalization of the sketched vectors $SV_m$ and $S(AV_m-V_mT_m)$ (see line 3 of Algorithm~\ref{alg:rstgrksm}) is much cheaper than the orthogonalization of the original unsketched vectors. Specifically, with $\epsilon = 10^{-0.5}$ so that the sketched subspace dimension is $s = 10m$, sketching would bring no benefit at all if the problem size is merely $n=10000$ and the maximum subspace dimension $m$ is as large as $1000$, since $s=10m=n$ here. Also, the saving in the cost of orthogonalization would lead to a considerable reduction of runtime if the linear solves needed at each RKSM step are relatively cheap. {Moreover, though the cost of selective orthogonalization with our parameters below also increases quadratically with $m$, it does so with a much smaller prefactor than the standard orthogonalization, thanks to the inner products between sketched basis vectors replacing original vectors, orthogonalization against fewer basis vectors, and no need for reorthogonalization. A larger subspace dimension $m$ would highlight the advantage of the sketched method.}

For each test matrix, we chose four different blocks $B$. For the first and the third choices of $B$ (e.g., $B=\mathrm{randn}(n,15)$ and $B=A^{-1}\widehat{I}_{90}$ for problem \textit{gyro\_m}), the largest rational Krylov subspace dimension is set to $m=500$, whereas for the second and the fourth choices of $B$ ($B = \mathrm{randn}(n,30)$ and $B=A^{-1}\widehat{I}_{300}$ for \textit{gyro\_m}), we let $m = 1000$. The depth of selective orthogonalization of the sketched algorithm is $\frac{m}{4}=0.25m$ for most problems, except for problems \textit{aerofoil} with $B = A^{-1}\widehat{I}_{60}$, and \textit{chipcool1} with $B = A^{-1}\widehat{I}_{40}$ and $B = A^{-1}\widehat{I}_{80}$, where the orthogonalization depth is $0.35m$.

Table~\ref{test3} shows that problems \textit{rail\_79841}, \textit{wing}, \textit{delaunay\_n16}, and \textit{Laplace2D} satisfy the features mentioned above for the sketched algorithm to outperform significantly in runtime. Problems \textit{filter3D}, \textit{boneS01}, and \textit{Laplace3D} have a medium-high cost in linear solves, but they are large in problem size; when a very large subspace is needed ($m \approx 900$--$1000$), the sketched algorithm also has a material advantage. For these problems with $B$ of higher rank, the runtime of the sketched algorithm is about $13$--$35\%$ of that of the standard algorithm. In particular, for \textit{Laplace2D} with $B = \mathrm{randn}(n,60)$, the runtime of the sketched method is only one eighth of that for the standard method. Moreover, the runtime ratio between the sketched and the standard algorithm tends to decrease whenever the largest subspace dimension $m$ increases. Meanwhile, problems \textit{chipcool1} and \textit{aerofoil} are small and costly in linear solves; they do not satisfy any features that would benefit the sketched method. Similarly, \textit{gas\_sensor} is of medium-large size, but the linear solves are very expensive. For these problems, the saving in runtime is modest, limited, or even nonexistent (e.g. \textit{gas\_sensor} with $B = C(:,1:14)$) for the sketched method. In general, the sketched method usually outperforms in runtime for large problems.

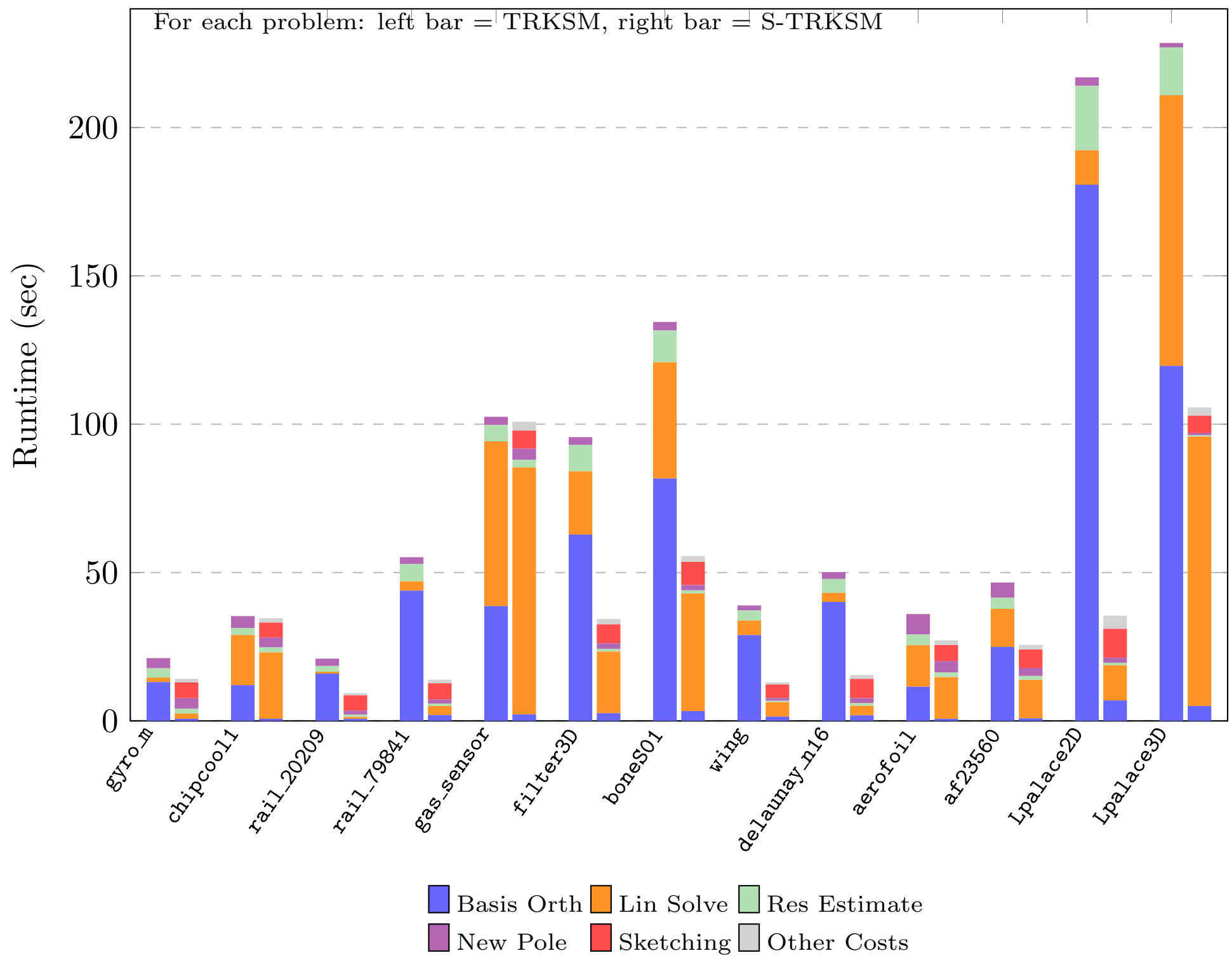
\begin{figure}[htbp!]
\centering
\begin{tikzpicture}

\def\cOrth{blue!60}
\def\cSolve{orange!85}
\def\cRes{green!60!black!30}
\def\cPole{violet!60}
\def\cSketch{red!70}
\def\cOther{gray!35}

\pgfplotsset{
  runtimeplot/.style={
    ybar stacked,
    bar width=7pt,
    enlarge x limits=0,
    width=\linewidth,
    height=9cm,
    ymin=0,
    ymax=240,
    enlarge y limits=false,
    ylabel={Runtime (sec)},
    ymajorgrids=true,
    grid style={dashed},
    legend style={font=\scriptsize, at={(0.5,-0.22)}, anchor=north, legend columns=3, draw=none, legend cell align={left}},
    symbolic x coords={
      LeftGap,
      gyroStd,gyroSk,gyroGap,
      chipStd,chipSk,chipGap,
      rail20209Std,rail20209Sk,rail20209Gap,
      rail79841Std,rail79841Sk,rail79841Gap,
      gasStd,gasSk,gasGap,
      filter3DStd,filter3DSk,filter3DGap,
      boneStd,boneSk,boneGap,
      wingStd,wingSk,wingGap,
      delStd,delSk,delGap,
      aeroStd,aeroSk,aeroGap,
      afStd,afSk,afGap,
      L2Std,L2Sk,L2Gap,
      L3Std,L3Sk,L3Gap
    },
    xtick={
      gyroStd, chipStd, rail20209Std, rail79841Std, gasStd, filter3DStd, boneStd,
      wingStd, delStd, aeroStd, afStd, L2Std, L3Std
    },
    xticklabels={
      \texttt{gyro\_m}, \texttt{chipcool1}, \texttt{rail\_20209}, \texttt{rail\_79841},
      \texttt{gas\_sensor}, \texttt{filter3D}, \texttt{boneS01}, \texttt{wing},
      \texttt{delaunay\_n16}, \texttt{aerofoil}, \texttt{af23560}, \texttt{Lpalace2D}, \texttt{Lpalace3D}
    },
    xticklabel style={font=\scriptsize, rotate=55, anchor=east},
  }
}

\begin{axis}[runtimeplot]

\addplot[fill=\cOrth, draw=none] coordinates {
(LeftGap,0)
(gyroStd,13.086) (gyroSk,0.714) (gyroGap,0)
(chipStd,12.034) (chipSk,0.753) (chipGap,0)
(rail20209Std,15.929) (rail20209Sk,0.693) (rail20209Gap,0)
(rail79841Std,43.895) (rail79841Sk,1.956) (rail79841Gap,0)
(gasStd,38.706) (gasSk,2.175) (gasGap,0)
(filter3DStd,62.790) (filter3DSk,2.590) (filter3DGap,0)
(boneStd,81.687) (boneSk,3.305) (boneGap,0)
(wingStd,28.914) (wingSk,1.438) (wingGap,0)
(delStd,40.089) (delSk,1.860) (delGap,0)
(aeroStd,11.505) (aeroSk,0.697) (aeroGap,0)
(afStd,24.928) (afSk,0.876) (afGap,0)
(L2Std,180.704) (L2Sk,6.900) (L2Gap,0)
(L3Std,119.653) (L3Sk,5.026) (L3Gap,0)
};
\addlegendentry{Basis Orth}

\addplot[fill=\cSolve, draw=none] coordinates {
(LeftGap,0)
(gyroStd,1.503) (gyroSk,1.743) (gyroGap,0)
(chipStd,16.908) (chipSk,22.319) (chipGap,0)
(rail20209Std,0.638) (rail20209Sk,0.632) (rail20209Gap,0)
(rail79841Std,3.156) (rail79841Sk,3.052) (rail79841Gap,0)
(gasStd,55.514) (gasSk,83.247) (gasGap,0)
(filter3DStd,21.315) (filter3DSk,20.765) (filter3DGap,0)
(boneStd,39.189) (boneSk,39.636) (boneGap,0)
(wingStd,4.890) (wingSk,4.828) (wingGap,0)
(delStd,3.069) (delSk,3.206) (delGap,0)
(aeroStd,13.997) (aeroSk,14.018) (aeroGap,0)
(afStd,12.833) (afSk,12.923) (afGap,0)
(L2Std,11.599) (L2Sk,11.848) (L2Gap,0)
(L3Std,91.232) (L3Sk,90.777) (L3Gap,0)
};
\addlegendentry{Lin Solve}

\addplot[fill=\cRes, draw=none] coordinates {
(LeftGap,0)
(gyroStd,3.124) (gyroSk,1.704) (gyroGap,0)
(chipStd,2.397) (chipSk,1.710) (chipGap,0)
(rail20209Std,1.964) (rail20209Sk,0.770) (rail20209Gap,0)
(rail79841Std,5.802) (rail79841Sk,0.820) (rail79841Gap,0)
(gasStd,5.505) (gasSk,2.559) (gasGap,0)
(filter3DStd,8.906) (filter3DSk,0.934) (filter3DGap,0)
(boneStd,10.703) (boneSk,1.074) (boneGap,0)
(wingStd,3.417) (wingSk,0.524) (wingGap,0)
(delStd,4.681) (delSk,0.897) (delGap,0)
(aeroStd,3.632) (aeroSk,1.651) (aeroGap,0)
(afStd,3.734) (afSk,1.329) (afGap,0)
(L2Std,21.736) (L2Sk,0.838) (L2Gap,0)
(L3Std,16.086) (L3Sk,0.464) (L3Gap,0)
};
\addlegendentry{Res Estimate}

\addplot[fill=\cPole, draw=none] coordinates {
(LeftGap,0)
(gyroStd,3.441) (gyroSk,3.495) (gyroGap,0)
(chipStd,3.977) (chipSk,3.288) (chipGap,0)
(rail20209Std,2.449) (rail20209Sk,1.359) (rail20209Gap,0)
(rail79841Std,2.302) (rail79841Sk,1.386) (rail79841Gap,0)
(gasStd,2.772) (gasSk,3.788) (gasGap,0)
(filter3DStd,2.620) (filter3DSk,1.847) (filter3DGap,0)
(boneStd,2.858) (boneSk,1.779) (boneGap,0)
(wingStd,1.687) (wingSk,0.943) (wingGap,0)
(delStd,2.291) (delSk,1.736) (delGap,0)
(aeroStd,6.872) (aeroSk,3.799) (aeroGap,0)
(afStd,5.143) (afSk,2.717) (afGap,0)
(L2Std,2.836) (L2Sk,1.725) (L2Gap,0)
(L3Std,1.479) (L3Sk,0.779) (L3Gap,0)
};
\addlegendentry{New Pole}

\addplot[fill=\cSketch, draw=none] coordinates {
(LeftGap,0)
(gyroStd,0) (gyroSk,5.297) (gyroGap,0)
(chipStd,0) (chipSk,5.064) (chipGap,0)
(rail20209Std,0) (rail20209Sk,5.156) (rail20209Gap,0)
(rail79841Std,0) (rail79841Sk,5.469) (rail79841Gap,0)
(gasStd,0) (gasSk,6.080) (gasGap,0)
(filter3DStd,0) (filter3DSk,6.406) (filter3DGap,0)
(boneStd,0) (boneSk,7.810) (boneGap,0)
(wingStd,0) (wingSk,4.517) (wingGap,0)
(delStd,0) (delSk,6.451) (delGap,0)
(aeroStd,0) (aeroSk,5.401) (aeroGap,0)
(afStd,0) (afSk,6.205) (afGap,0)
(L2Std,0) (L2Sk,9.720) (L2Gap,0)
(L3Std,0) (L3Sk,5.814) (L3Gap,0)
};
\addlegendentry{Sketching}

\addplot[fill=\cOther, draw=none] coordinates {
(LeftGap,0)
(gyroStd,0.022) (gyroSk,1.248) (gyroGap,0)
(chipStd,0.016) (chipSk,1.498) (chipGap,0)
(rail20209Std,0.011) (rail20209Sk,0.705) (rail20209Gap,0)
(rail79841Std,0.033) (rail79841Sk,1.199) (rail79841Gap,0)
(gasStd,0.033) (gasSk,2.959) (gasGap,0)
(filter3DStd,0.048) (filter3DSk,1.822) (filter3DGap,0)
(boneStd,0.084) (boneSk,1.965) (boneGap,0)
(wingStd,0.027) (wingSk,0.735) (wingGap,0)
(delStd,0.031) (delSk,1.321) (delGap,0)
(aeroStd,0.031) (aeroSk,1.592) (aeroGap,0)
(afStd,0.016) (afSk,1.587) (afGap,0)
(L2Std,0.091) (L2Sk,4.418) (L2Gap,0)
(L3Std,0.098) (L3Sk,2.794) (L3Gap,0)
};
\addlegendentry{Other Costs}

\node[font=\scriptsize, anchor=west] at (rel axis cs:0.01,0.98)
{For each problem: left bar = TRKSM, right bar = S-TRKSM};

\end{axis}
\end{tikzpicture}
\caption{Runtime breakdown (stacked) for TRKSM (left bar) and S-TRKSM (right bar) on 13 problems, with matrix $B$ as the first choice of the four  in Table~\ref{test3} leading to convergence ($B=\mathrm{randn(n,p)}$ for most problems). Other cost accounts for unreported overhead so each stack matches the total runtime.}
\label{fig:all_runtime_stacked}\vspace{-0.1in}
\end{figure}

\begin{figure}[hbtp!]
\centering
\setlength{\tabcolsep}{2pt} %
\renewcommand{\arraystretch}{0} %

\begin{tabular}{@{}c c@{}}
  \includegraphics[width=.50\textwidth]{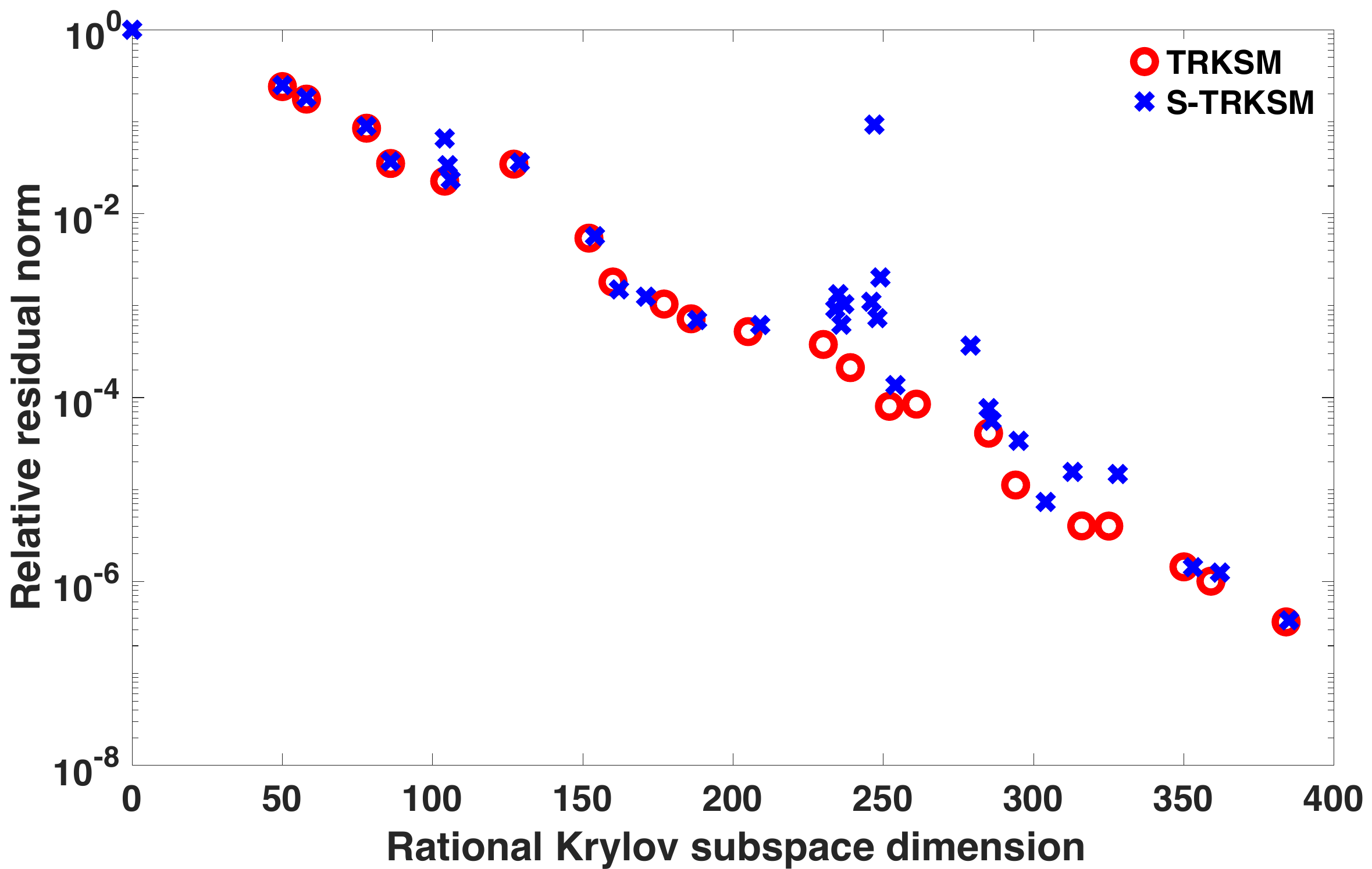} &
  \includegraphics[width=.49\textwidth]{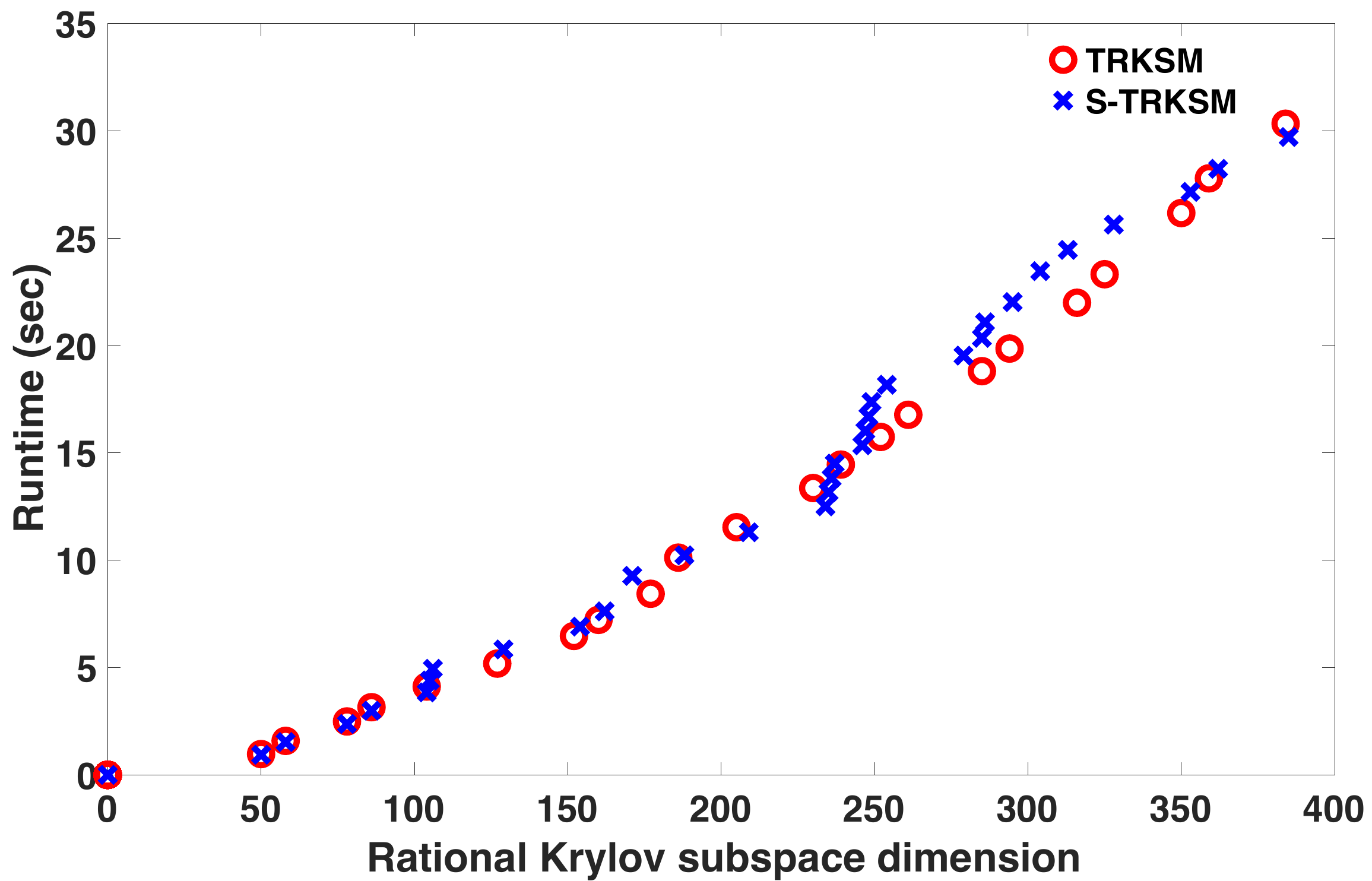} \\
  \includegraphics[width=.49\textwidth]{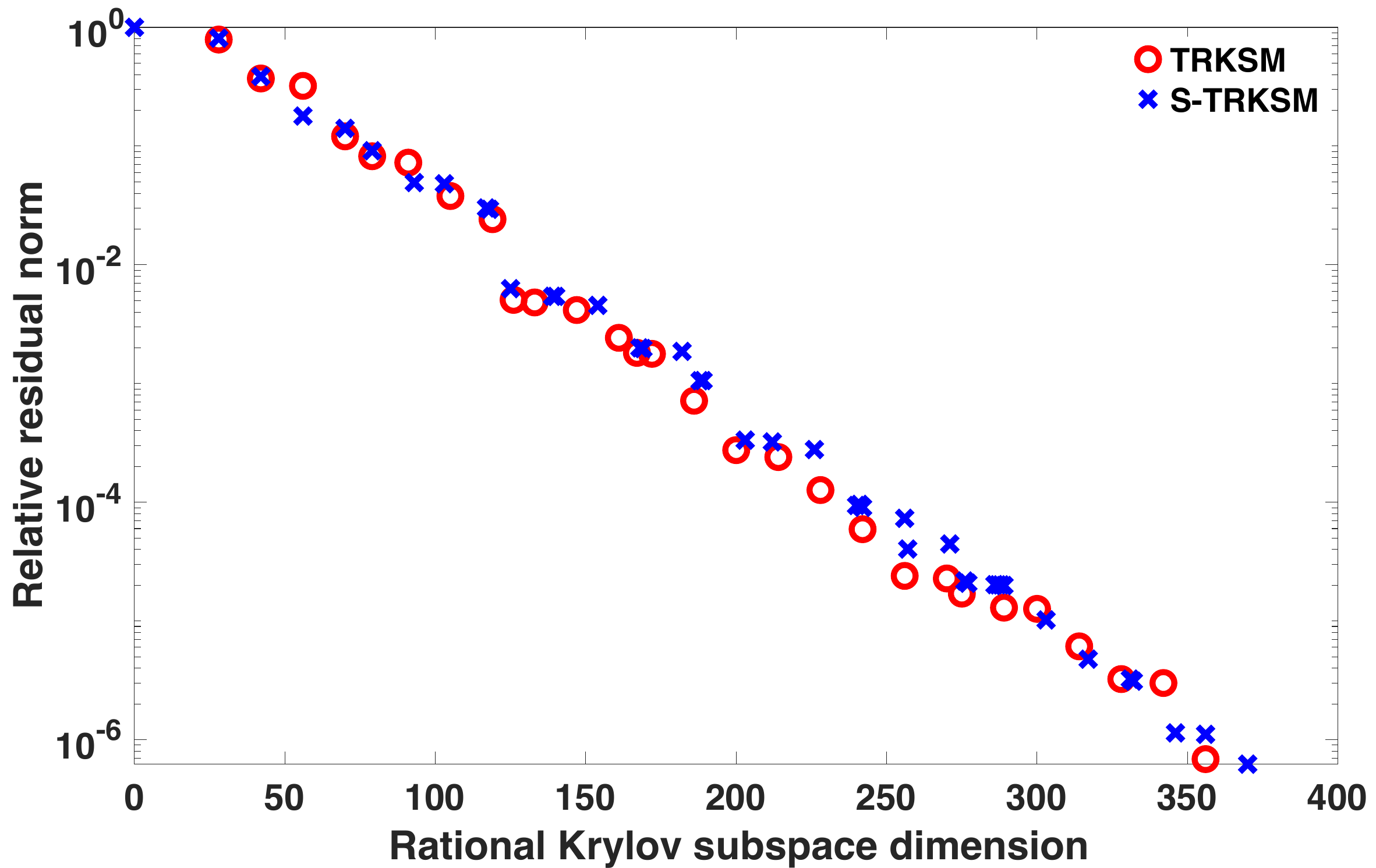} &
  \includegraphics[width=.49\textwidth]{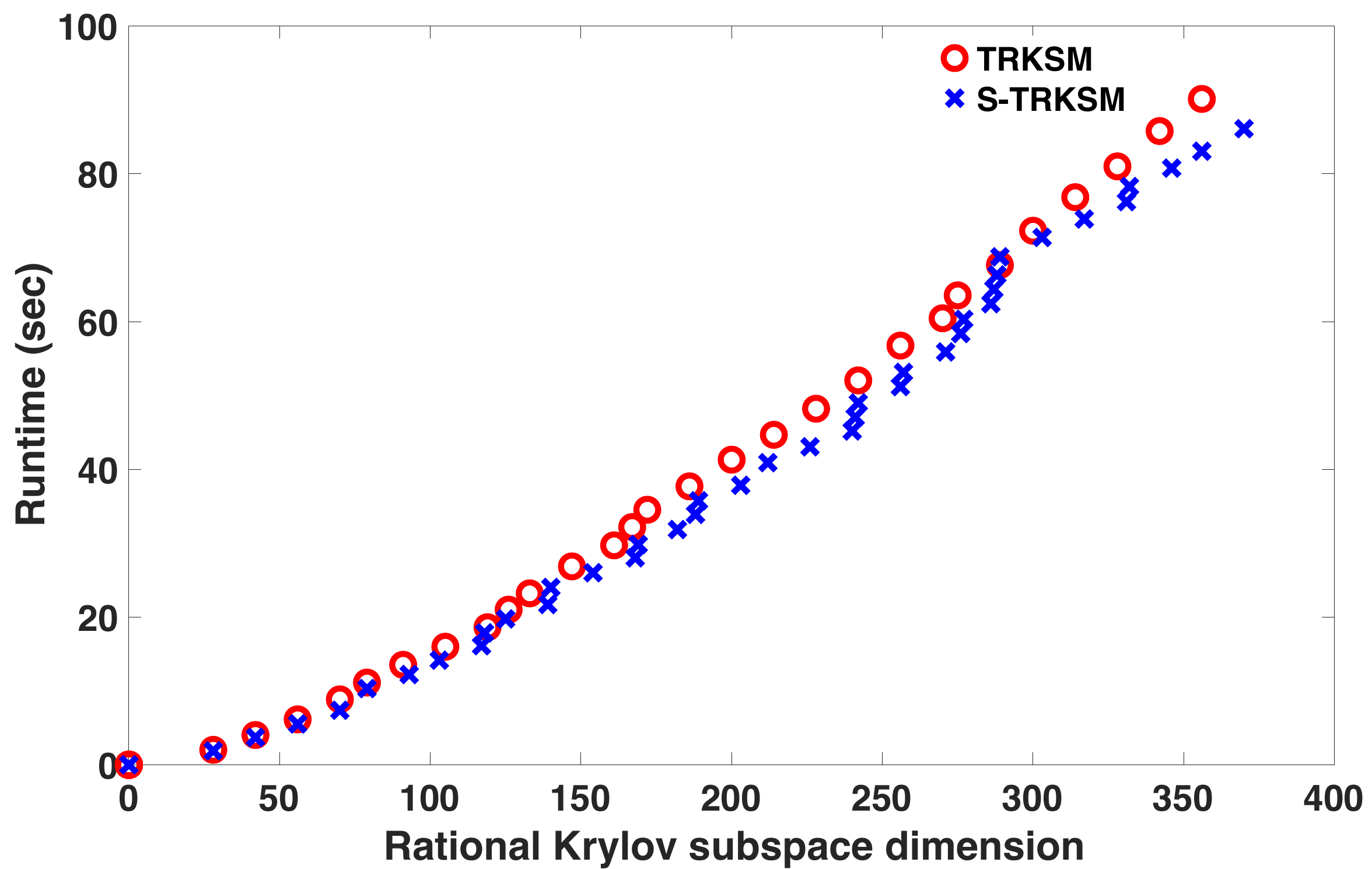}
\end{tabular}

\caption{{Equation relative residual norm $\|AX_m + X_m A^T+ BB^T\|_F/\|BB^T\|_F$ and runtime (secs) versus rational Krylov subspace dimension. Top\,:\,problem \textit{chipcool1} with $B=[C\:\widehat{E}_{20}]$; bottom\,:\,problem \textit{gas\_sensor} with $B=C(:,1\!:\!14)$} (see Table~\ref{test3})}
\label{fig:restimevsdim}\vspace{-0.2in}
\end{figure}

{In Figure~\ref{fig:all_runtime_stacked}, a stacked bar plot is drawn to show the details of the runtime of different computational tasks. For each problem, we let $B$ be the first of the four choices in Table~\ref{test3} leading to convergence of both methods. It is evident that the sketched algorithm significantly reduced the runtime of basis orthogonalization for all problems. In addition, our fast equation norm estimate discussed in section~\ref{subsect:resest} is also substantially faster than the counterpart of the standard algorithm. For problems \textit{chipcool1} and \textit{gas\_sensor}, the total runtime is about the same for both algorithms because the sketched method takes considerably more iteration steps and hence a longer time on the linear solves with distinct poles. The cost of sketching is modest compared to the total runtime, especially for larger problems where the sketched vectors have much fewer elements than the original vectors ($s \ll n$).}

There are two problems for which the sketched method needs considerably more iteration steps, i.e., more linear solves with distinct poles, even though the final subspace dimensions are comparable. {For \textit{chipcool1} with \(B=[C\ \widehat{E}_{20}]\), the standard algorithm takes $26$ steps whereas the sketched variant takes $38$ steps. The same effect is even more visible for \textit{gas\_sensor}: with \(B=C(:,1\!:\!14)\), the number of steps increases from $31$ to $49$, and with \(B=C\), from $32$ to $51$. Since \textit{gas\_sensor} also has a high sparse-factor density, the additional distinct-pole linear solves dominate the runtime and offset the savings from sketched orthogonalization and residual estimation. This explains why the linear-solve portion of the sketched runtime is larger for \textit{gas\_sensor} in Figure~\ref{fig:restimevsdim}. We speculate that the universal threshold used to select effective tangential directions to expand the rational Krylov subspace could be lowered for the sketched method, at least for these two problems, so that it may expand the subspace with more vectors each step and save the number of iteration steps and runtime. Fortunately, for most other high-rank examples the iteration counts of the two algorithms are similar, and the reductions in orthogonalization and residual-estimation costs translate into a clear runtime advantage.} In addition, there are a few problems for which the standard algorithm was not able to satisfy the prescribed tolerances and is marked as `FAIL' in Table~\ref{test3}, but the sketched algorithm achieved the tolerances without difficulty.

\section{Conclusion}
We proposed a randomly sketched tangential rational Krylov subspace method for solving large continuous-time Lyapunov equations $AX+XA^T+BB^T=0$ with $\mathrm{rank}(B)\gg 1$. To the best of our knowledge, this is the first RKSM variant that uses a non-orthonormal basis supported by randomized sketching. Sketching reduces the cost of pole selection, residual-norm estimation, and basis orthogonalization. We established an exactness property in the sketched case analogous to that of the standard (non-sketched) RKSM and derived a new error bound that makes the dependence on the poles explicit. This bound supports linear (exponential) convergence and is efficiently computable via the fast Fourier transform (FFT).

Numerical results (Tables~\ref{test2} and \ref{test3}) indicate that the sketched method can deliver more accurate approximations with smaller residual norms than the standard variant. For Lyapunov equations with low-rank $B$, using a few repeated poles with sparse factorizations can be competitive in runtime. When $B$ has higher rank and RKSM requires a larger subspace, the sketched method with adaptive poles is faster, especially for large problems where linear solves are relatively cheap. Provided storage is affordable, it can therefore solve high-rank-$B$ Lyapunov equations more efficiently.

This work contributes to recent efforts to accelerate long-recurrence iterative methods using randomized sketching. The proposed algorithm enables reduced computation time for challenging problems in hydrodynamic stability, model reduction, and control with many inputs and outputs. \medskip

{\bf Acknowledgement.} We thank the referees for their careful reading of our manuscript and making constructive comments. GPT-5.2 was used to polish our written text for spelling, grammar, and general style, and to generate Figures \ref{fig:runtime_stacked_5problems} and \ref{fig:all_runtime_stacked}. 

\bibliographystyle{siam}
\bibliography{ref}

\end{document}